# ON THE MAXIMUM QUEUE LENGTH IN THE SUPERMARKET MODEL

By Malwina J. Luczak and Colin McDiarmid

*London School of Economics and University of Oxford*

There are $n$ queues, each with a single server. Customers arrive in a Poisson process at rate $\lambda n$, where $0 < \lambda < 1$. Upon arrival each customer selects $d \geq 2$ servers uniformly at random, and joins the queue at a least-loaded server among those chosen. Service times are independent exponentially distributed random variables with mean 1. We show that the system is rapidly mixing, and then investigate the maximum length of a queue in the equilibrium distribution. We prove that with probability tending to 1 as $n \to \infty$ the maximum queue length takes at most two values, which are $\ln \ln n / \ln d + O(1)$.

**1. Introduction.** We study a well-known queueing model with $n$ separate queues, each with a single server. Customers arrive into the system in a Poisson process at rate $\lambda n$, where $0 < \lambda < 1$ is a constant. Upon arrival each customer chooses $d$ queues uniformly at random with replacement, and joins a shortest queue amongst those chosen (where she breaks ties by choosing the first of the shortest queues in the list of $d$). Here $d$ is a fixed positive integer. Customers are served according to the first-come–first-served discipline. Service times are independent exponentially distributed random variables with mean 1.

A number of authors have studied this model before, as well as its extension to a Jackson network setting [2, 3, 4, 11, 13, 14, 15, 17]. For instance, it is shown in [2] that the system is *chaotic*, provided that it starts close to a suitable deterministic initial state, or is in equilibrium. This means that the paths of members of any fixed finite subset of queues are asymptotically independent of one another, uniformly on bounded time intervals. This result implies a law of large numbers for the time evolution of the proportion of queues of different lengths, that is, for the empirical measure on the path space [2]. In particular for each fixed positive integer $k_0$, as $n$ tends to infinity the proportion of queues with length at least $k_0$ converges weakly (when the infinite-dimensional state space is endowed with the product topology) to a function $v_t(k_0)$, where $v_t(0) = 1$ for all $t \geq 0$ and $(v_t(k): k \in \mathbb{N})$ is the unique solution to the system of differential equations

$$\text{(1)} \qquad \frac{dv_t(k)}{dt} = \lambda(v_t(k-1)^d - v_t(k)^d) - (v_t(k) - v_t(k+1))$$

for $k \in \mathbb{N}$. Here we assume appropriate initial conditions $(v_0(k): k \in \mathbb{N})$ such that $1 \geq v_0(1) \geq v_0(2) \geq \cdots \geq 0$. Further, again for a fixed positive integer $k_0$, as $n$ tends to infinity, in the equilibrium distribution this proportion converges in probability to $\lambda^{1+d+\cdots+d^{k_0-1}}$, and thus the probability that a given queue has length at least $k_0$ also converges to $\lambda^{1+d+\cdots+d^{k_0-1}}$.

Although these results refer only to fixed queue length $k_0$ and bounded time intervals, they suggest that when $d \geq 2$, in equilibrium the maximum queue length may usually be $O(\ln \ln n)$. Our main contribution is to show that this is indeed the case, and to give precise results on the behavior of the maximum queue length. In particular, we shall see that when $d \geq 2$, with probability tending to 1 as $n \to \infty$, in the equilibrium distribution the maximum queue length takes at most two values; and these values are $\ln \ln n / \ln d + O(1)$. We show also that the system is rapidly mixing, that is, the distribution settles down quickly to the equilibrium distribution. Another natural question concerns fluctuations when in the equilibrium distribution: how long does it take to see large deviations of the maximum queue length from its stationary median? We provide an answer by establishing strong concentration estimates over time intervals of length polynomial in $n$. Our techniques are partly combinatorial, and are used also in [7, 8, 9]. In particular, in [8] we use the concentration estimates obtained here to establish quantitative results on the convergence of the distribution of a queue length and on the asymptotic independence of small subsets of queues, the "chaotic behavior" of the system.

Recently, in [6, 10], a quantitative approximation has been obtained for the supermarket model, including a law of large numbers and a central limit theorem. These results rely on properties of continuous-time exponential martingales and strong approximation of Poisson processes by Brownian motion. A "localization" technique yields tight bounds on the deviation probabilities uniformly in all co-ordinates of the infinite-dimensional state space, in a spirit somewhat akin to the approach adopted in this paper and in [7, 8, 9]. However, the results in [6, 10] concern solely fixed-length time intervals and do not extend to equilibrium behavior.

Let us introduce some notation so that we can state our results. Consider the $n$-queue model. For each time $t \geq 0$ and each $j = 1, \ldots, n$ let $X_t^{(n)}(j)$ or $X_t(j)$ denote the number of customers in queue $j$, always including the customer currently being served if there is one. (We shall keep the superscript "$n$" in the notation in this section, but then drop it in later sections.) We make the usual assumption of right-continuity of the sample paths. Let $X_t^{(n)}$ or $X_t$ be the *queue-lengths vector* $(X_t^{(n)}(1), \ldots, X_t^{(n)}(n))$. For a given positive integer $n$, $X_t^{(n)}$ is an ergodic continuous-time Markov chain. Thus there is a unique stationary distribution $\mathbf{\Pi}^{(n)}$ or $\mathbf{\Pi}$ for the vector of queue lengths; and, whatever the distribution of the starting state, the distribution of the queue-lengths vector $X_t^{(n)}$ at time $t$ converges to $\mathbf{\Pi}^{(n)}$ as $t \to \infty$. We will show that, with reasonable initial conditions, the convergence is very fast. Note that the $L_1$-norm $\|X_t\|_1$ of $X_t$ is the total number of customers present at time $t$, and the $L_\infty$-norm $\|X_t\|_\infty$ is the maximum queue length.

The probability law or distribution of a random variable $X$ will be denoted by $\mathcal{L}(X)$. The *total variation distance* between two probability distributions $\mu_1$ and $\mu_2$ may be defined by $d_{\text{TV}}(\mu_1, \mu_2) = \sup_A |\mathbf{Pr}(X \in A) - \mathbf{Pr}(Y \in A)|$, where the supremum is over all events $A$, or equivalently by $d_{\text{TV}}(\mu_1, \mu_2) = \inf \mathbf{Pr}(X \neq Y)$, where the infimum is over all couplings of $X$ and $Y$ where $\mathcal{L}(X) = \mu_1$ and $\mathcal{L}(Y) = \mu_2$.

For any given state $x$ we shall write $\mathcal{L}(X_t^{(n)}, x)$ to denote the law of $X_t^{(n)}$ given $X_0^{(n)} = x$. For $\varepsilon > 0$, the *mixing time* $\tau^{(n)}(\varepsilon, x)$ starting from state $x$ is defined by

$$\tau^{(n)}(\varepsilon, x) = \inf\{t \geq 0 : d_{\text{TV}}(\mathcal{L}(X_t^{(n)}, x), \mathbf{\Pi}^{(n)}) \leq \varepsilon\}.$$

Our first theorem shows that if we start from any initial state in which the queues are not too long, then the mixing time is small. In particular, if $\varepsilon > 0$ is fixed and $\mathbf{0}$ denotes the all-zero $n$-vector, then $\tau^{(n)}(\varepsilon, \mathbf{0})$ is $O(\ln n)$. More generally, this holds if $\varepsilon > 0$ is not too small (explicitly, if $\varepsilon^{-1}$ is bounded polynomially in $n$), $\|x\|_1$ is $O(n)$ and $\|x\|_\infty$ is $O(\ln n)$. Observe that the quantity $\delta_{n,t}$ below is 0 if the initial state is $\mathbf{0}$.

THEOREM 1.1.  *Let $0 < \lambda < 1$ and let $d$ be a fixed positive integer. For each constant $c > 0$ there exists a constant $\eta > 0$ such that the following holds for each positive integer $n$. Consider any distribution of the initial queue-lengths vector $X_0^{(n)}$, and for each time $t \geq 0$ let*

$$\delta_{n,t} = \mathbf{Pr}(\|X_0^{(n)}\|_1 > cn) + \mathbf{Pr}(\|X_0^{(n)}\|_\infty > \eta t).$$

*Then*

$$d_{\text{TV}}(\mathcal{L}(X_t^{(n)}), \mathbf{\Pi}^{(n)}) \leq ne^{-\eta t} + 2e^{-\eta n} + \delta_{n,t}.$$

The $O(\ln n)$ upper bound on the mixing time $\tau$ is of the right order. Indeed, we shall see that for a suitable constant $\theta > 0$, if $t \leq \theta \ln n$, then

$$(2) \qquad d_{\text{TV}}(\mathcal{L}(X_t^{(n)}), \mathbf{\Pi}^{(n)}) = 1 - e^{-\Omega(\ln^2 n)}.$$

Thus $\tau^{(n)}(\varepsilon, \mathbf{0})$ is $\Theta(\ln n)$ as long as both $\varepsilon^{-1}$ and $(1-\varepsilon)^{-1}$ are bounded polynomially in $n$.

Our primary interest is in the maximum queue length $M_t^{(n)} = \|X_t^{(n)}\|_\infty$. Since the system mixes rapidly it is natural to consider the stationary case. Our model exhibits the "power of two choices" phenomenon (see, e.g., [15]); that is, when we move from $d = 1$ choice to $d = 2$ choices, the typical maximum queue length $M_t^{(n)}$ drops dramatically.

We are most interested in the case $d \geq 2$, but first we consider the easy case $d = 1$ in order to set the scene. This case is straightforward, since in equilibrium the queue lengths are i.i.d. geometric random variables with parameter $\lambda$. We find that the maximum queue length $M_t^{(n)}$ is about $\frac{\ln n}{\ln(1/\lambda)}$, and it is not concentrated on a bounded range of values (in contrast to the behavior in the balls-and-bins model [7], where the maximum load is concentrated on at most two adjacent values, even when $d = 1$).

Given a sequence of events $A_1, A_2, \ldots$, we say that $A_n$ holds *asymptotically almost surely* (a.a.s.) if $A_n$ holds with probability tending to 1 as $n \to \infty$.

THEOREM 1.2.  *Let $0 < \lambda < 1$ and let $d = 1$. For each positive integer $n$, suppose that the queue-lengths vector $X_0^{(n)}$ is in the stationary distribution (and thus so is the maximum queue length $M_t^{(n)}$).*

(a) *For each nonnegative integer $m$*

$$\mathbf{Pr}(M_t^{(n)} \leq m) = (1 - \lambda^{m+1})^n.$$

*Thus if $m = m(n)$ and $n \to \infty$, then $M_t^{(n)} \geq m(n)$ a.a.s. if and only if $m(n) - \frac{\ln n}{\ln(1/\lambda)} \to -\infty$; and $M_t^{(n)} \leq m(n)$ a.a.s. if and only if $m(n) - \frac{\ln n}{\ln(1/\lambda)} \to +\infty$.*

(b) *For any subexponential time $\tau \geq 0$ (i.e., $\tau = e^{n^{o(1)}}$),*

$$\left(\min_{0 \leq t \leq \tau} M_t^{(n)}\right) \frac{\ln(1/\lambda)}{\ln n} \to 1 \qquad \text{in probability as } n \to \infty.$$

(c) *For any constant $K > 0$,*

$$\left(\max_{0 \leq t \leq n^K} M_t^{(n)}\right) \frac{\ln(1/\lambda)}{\ln n} \to K + 1 \qquad \text{in probability as } n \to \infty.$$

Now we consider the case $d \geq 2$, when the maximum queue length $M_t^{(n)}$ is far smaller, and it is concentrated on just two values $m_d$ and $m_d - 1$. This is our main result.

THEOREM 1.3. *Let $0 < \lambda < 1$ and let $d \geq 2$ be an integer. Then there exists an integer-valued function $m_d = m_d(n) = \ln \ln n / \ln d + O(1)$ such that the following holds. For each positive integer $n$, suppose that the queue-lengths vector $X_0^{(n)}$ is in the stationary distribution (and thus so is the maximum queue length $M_t^{(n)}$). Then for each time $t \geq 0$, $M_t^{(n)}$ is $m_d(n)$ or $m_d(n) - 1$ a.a.s.; and further, for any constant $K > 0$ there exists $c = c(K)$ such that*

$$(3) \qquad \max_{0 \leq t \leq n^K} |M_t^{(n)} - \ln \ln n / \ln d| \leq c \qquad a.a.s.$$

The functions $m_2(n), m_3(n), \ldots$ may be defined as follows. For $d = 2, 3, \ldots$ let $i_d(n)$ be the least integer $i$ such that $\lambda^{(d^i-1)/(d-1)} < n^{-1/2} \ln^2 n$. Then we let $m_2(n) = i_2(n) + 1$, and for $d \geq 3$ let $m_d(n) = i_d(n)$. (We shall see that with high probability the proportion of queues of length at least $i$ is close to $\lambda^{(d^i-1)/(d-1)}$.)

Results can be obtained on deviations of the maximum queue length over longer time intervals, using arguments just as in [7]; we shall not discuss such results further here.

We have now described our model and stated our main results. The later sections of this paper are organized as follows. In Section 2 we prove Theorem 1.1, which shows that for each fixed positive integer $d$ the process $(X_t^{(n)})$ starting from a nice initial state mixes in logarithmic time. The proof is based on considering a pair of "adjacent" initial states and analyzing a suitable random walk. We give also a similar result involving the Wasserstein distance. In Section 3 we focus on the straightforward case $d = 1$, and prove Theorem 1.2. In Section 4 we show that a Lipschitz function of the queue-lengths vector in equilibrium is tightly concentrated around its mean. To do this, we consider a queue-lengths process $(X_t^{(n)})$ starting from **0**, and use the bounded differences approach to establish concentration at a suitable time $t$ when $\mathcal{L}(X_t^{(n)})$ is close to the equilibrium distribution. In Section 5 we use the concentration property to estimate the proportions of queues of at least some given lengths, and to bound their fluctuations over long time intervals. In the short Section 6 that follows we establish the logarithmic lower bound (2) on the mixing times in Theorem 1.1. In Section 7 we prove Theorem 1.3, and thus complete the proofs of the new results stated above. Finally, we make some brief concluding remarks.

Several times we shall use the fact that, if $X$ is a binomial or Poisson random variable with mean $\mu$, then for each $0 \leq \varepsilon \leq 1$ we have

$$(4) \qquad \mathbf{Pr}(X - \mu \leq -\varepsilon \mu) \leq e^{-(1/2)\varepsilon^2 \mu}$$

and

$$(5) \qquad \mathbf{Pr}(X - \mu \geq \varepsilon \mu) \leq e^{-(1/3)\varepsilon^2 \mu};$$

and if $x \geq 2e\mu$, then

(6) $$\mathbf{Pr}(X \geq x) \leq 2^{-x}.$$

For the above results, see, for example, Theorem 2.3(b) and inequalities (2.7) and (2.8) in [12].

**2. Rapid mixing: proof of Theorem 1.1.** In this section we shall in fact prove both Theorem 1.1 and a similar result involving the Wasserstein distance instead of the total variation distance. The *Wasserstein distance* may be defined by $d_\mathrm{W}(\mu_1, \mu_2) = \inf \mathbf{E}[\|X - Y\|_1]$ where the infimum is over all couplings where $\mathcal{L}(X) = \mu_1$ and $\mathcal{L}(Y) = \mu_2$. Observe that for integer-valued random variables, the total variation distance between the corresponding laws is always at most the Wasserstein distance. The following result will also be used in [8], where we consider the asymptotic distribution of a queue length.

LEMMA 2.1. *Let $0 < \lambda < 1$ and let $d$ be a fixed positive integer. For each constant $c > \frac{\lambda}{1-\lambda}$ there exists a constant $\eta > 0$ such that the following holds for each positive integer $n$. Let $M$ denote the stationary maximum queue length. Consider any distribution of the initial queue-lengths vector $X_0$ such that $\|X_0\|_1$ has finite mean. For each time $t \geq 0$ let*

$$\delta_{n,t} = 2\mathbf{E}[\|X_0\|_1 \mathbf{1}_{\|X_0\|_1 > cn}] + 2cn\mathbf{Pr}(M_0 > \eta t).$$

*Then*

$$d_\mathrm{W}(\mathcal{L}(X_t), \mathbf{\Pi}) \leq ne^{-\eta t} + 2cn\mathbf{Pr}(M > \eta t) + 2e^{-\eta n} + \delta_{n,t}.$$

In order to prove this result we shall couple the queue-lengths process $(X_t)$ and a corresponding copy $(Y_t)$ of the process in equilibrium in such a way that with high probability $\|X_t - Y_t\|_1$ decreases quickly to 0. We model departures by a Poisson process at rate $n$ together with an independent selection process that generates a uniformly random queue at each event time. If the queue selected is nonempty, then the customer currently in service departs; otherwise nothing happens.

We start the proof with a lemma concerning the return times to the origin of a generalized random walk on $\{0, 1, 2, \ldots\}$. This lemma will be needed later to show that a certain coupling happens quickly.

LEMMA 2.2. *Let $\phi_0, \phi_1, \phi_2, \ldots$ be a filtration. Let $Z_1, Z_2, \ldots$ be $\{0, \pm 1\}$-valued random variables, where each $Z_i$ is $\phi_i$-measurable. Let $S_0 \geq 0$ a.s. and for each positive integer $j$ let $S_j = S_0 + \sum_{i=1}^{j} Z_i$. Let $A_0, A_1, \ldots$ be events, where each $A_i$ is $\phi_i$-measurable.*

*Suppose that there is a constant positive integer $k_0$ and a constant $\delta$ with $0 < \delta < 1/2$ such that*

$$\mathbf{Pr}(Z_i = -1|\phi_{i-1}) \geq \delta \qquad \text{on } A_{i-1} \cap \{S_{i-1} \in \{1, \ldots, k_0 - 1\}\}$$

*and*

$$\mathbf{Pr}(Z_i = -1|\phi_{i-1}) \geq \delta + 1/2 \qquad \text{on } A_{i-1} \cap \{S_{i-1} \geq k_0\}.$$

*Then there exists $\eta > 0$ such that for each positive integer $m$*

$$\mathbf{Pr}\left(\left(\bigcap_{i=1}^{m} S_i \neq 0\right) \cap \left(\bigcap_{i=0}^{m-1} A_i\right)\right) \leq \mathbf{Pr}(S_0 > \eta m) + 2e^{-\eta m}.$$

PROOF. Let us ignore the events $A_i$ in the meantime; we shall see later that it is easy to incorporate them into the argument. We define random times $I_j$ and $\hat{I}_j$ as follows. Let

$$I_0 = \min\{i \geq 1 : S_i < k_0\},$$

and let

$$I_{j+1} = \min\{i > I_j : S_i < k_0\}.$$

Further, let $\hat{I}_0 = I_0$ and let

$$\hat{I}_{j+1} = \min\{i \geq \hat{I}_j + k_0 : S_i < k_0\}.$$

Observe that always $\hat{I}_j \leq I_{k_0 j}$.

The key fact is that for all positive integers $m$ and $j$,

(7) $$\mathbf{Pr}\left(\bigcap_{i=1}^{m} S_i \neq 0\right) \leq (1 - \delta^{k_0})^j + \mathbf{Pr}(\hat{I}_j > m).$$

To see this, note first that if $\hat{I}_i = t$, then $S_t < k_0$ and, for $1 \leq k < k_0$, on $S_t = k$ we have

$$\mathbf{Pr}\left(\bigcap_{u=1}^{k}(Z_{t+u} = -1)|\phi_t\right) \geq \delta^k.$$

Now for each $i = 0, 1, \ldots$ let $B_i$ be the event that $S_r \neq 0$ for each $r = \hat{I}_i, \ldots, \hat{I}_i + k_0 - 1$. Then $\mathbf{Pr}(B_i|\phi_{\hat{I}_i}) \leq 1 - \delta^{k_0}$, and so $\mathbf{Pr}(\bigcap_{i=0}^{j-1} B_i) \leq (1 - \delta^{k_0})^j$. But

$$\mathbf{Pr}\left(\left(\bigcap_{i=1}^{m} S_i \neq 0\right) \cap (\hat{I}_j \leq m)\right) \leq \mathbf{Pr}\left(\bigcap_{i=0}^{j-1} B_i\right),$$

and (7) follows.

We need to investigate the times $\hat{I}_j$ in order to be able to ensure that the term $\mathbf{Pr}(\hat{I}_j > m)$ in (7) is small. First we consider the times $I_j$. Note that $\frac{r}{2} \leq (1-\delta)r(\frac{1}{2}+\delta)$ for $r \geq 0$. Let $h = \delta^2/4$. Then by (4), for all nonnegative integers $j$ and $r$,

$$\mathbf{Pr}(I_{j+1} - I_j > r | \phi_{I_j}) \leq \mathbf{Pr}\left(B\left(r, \frac{1}{2}+\delta\right) < \frac{r}{2}\right) \leq e^{-hr}.$$

(Here we are using $B$ to denote a binomial random variable.) Now let $s$ be a positive integer. Let $b = (\delta - 2\delta^2)^{-1}$, so $b > 0$. Note that for $r \geq bs$, we have $\frac{r+s}{2} \leq (1-\delta)r(\frac{1}{2}+\delta)$. Hence, again using (4), we see that, on the event $\{S_0 \leq s\}$,

$$\mathbf{Pr}(I_0 > r | \phi_0) \leq \mathbf{Pr}\left(B\left(r, \frac{1}{2}+\delta\right) < \frac{r+s}{2}\right) \leq e^{-hr}.$$

Thus in particular, for $r \geq bs$,

$$\mathbf{Pr}(I_0 \geq r) \leq e^{-hr} + \mathbf{Pr}(S_0 > s).$$

Let $Z$ be a random variable taking positive integer values and such that $\mathbf{Pr}(Z > r) = e^{-hr}$ for each positive integer $r$. Then $I_j - I_0$ is stochastically at most the sum of $j$ i.i.d. random variables, each distributed like $Z$. Note that $M_Z(h/2) = \mathbf{E}[e^{(h/2)Z}] < \infty$. Let $c_1 > 0$ be sufficiently large that $M_Z(h/2)e^{-c_1 h/2} < 1$, say $M_Z(h/2)e^{-c_1 h/2} = e^{-c_2}$ where $c_2 > 0$. Then

$$\mathbf{Pr}(I_j - I_0 > c_1 j) \leq e^{-(h/2)c_1 j} \mathbf{E}[e^{(h/2)(I_j - I_0)}]$$
$$\leq e^{-(h/2)c_1 j} M_Z(h/2)^j$$
$$= e^{-c_2 j}.$$

Putting the last two results together, and using the fact that $\hat{I}_j \leq I_{jk_0}$, we have for each $r \geq bs$,

$$\mathbf{Pr}(\hat{I}_j > r + c_1 j k_0) \leq \mathbf{Pr}(I_0 > r) + \mathbf{Pr}(I_{jk_0} - I_0 > c_1 j k_0)$$
$$\leq e^{-hr} + \mathbf{Pr}(S_0 > s) + e^{-c_2 j k_0}.$$

We may now complete the proof (still for the case without the events $A_i$). Let $\eta_1 > 0$ and $\eta_2 > 0$ be sufficiently small that $1 - \eta_1 - \eta_2 c_1 k_0 > 0$. Let $r = \lceil \eta_1 m \rceil$, $j = \lceil \eta_2 m \rceil$ and $s = \lfloor r/b \rfloor$. Then for $m$ sufficiently large, $r + c_1 j k_0 \leq m$; and so by (7) and the last inequality

$$\mathbf{Pr}\left(\bigcap_{i=1}^m S_i \neq 0\right) \leq e^{-\delta k_0 j} + e^{-hr} + \mathbf{Pr}(S_0 > s) + e^{-c_2 j k_0}$$
$$\leq e^{-\delta k_0 \eta_2 m} + e^{-h\eta_1 m} + \mathbf{Pr}(S_0 > s) + e^{-c_2 k_0 \eta_2 m}.$$

Thus there exists a constant $\eta' > 0$ and an integer $m_0$ such that for each integer $m \geq m_0$

$$\mathbf{Pr}\left(\bigcap_{i=1}^{m} S_i \neq 0\right) \leq \mathbf{Pr}(S_0 > \eta' m) + 2e^{-\eta' m}.$$

But now we may set $\eta = \min\{\eta', \ln 2/m_0\} > 0$ and then for each positive integer $m$

$$\mathbf{Pr}\left(\bigcap_{i=1}^{m} S_i \neq 0\right) \leq \mathbf{Pr}(S_0 > \eta m) + 2e^{-\eta m}.$$

Let us now bring in the events $A_i$. For each $i$ define $\tilde{Z}_i = Z_i \mathbb{I}_{A_i} - \mathbb{I}_{\overline{A_i}}$ (where $\overline{A_i}$ denotes the complement of $A_i$). Then $\mathbf{Pr}(\tilde{Z}_i = -1|\phi_{i-1}) \geq \delta$ on $\{S_{i-1} \in \{1,\ldots,k_0-1\}\}$ and $\mathbf{Pr}(\tilde{Z}_i = -1|\phi_{i-1}) \geq 1/2+\delta$ on $\{S_{i-1} \geq k_0\}$. Let $\tilde{S}_j = S_0 + \sum_{i=1}^{j} \tilde{Z}_i$. Then by what we have just proved applied to the $\tilde{Z}_i$,

$$\mathbf{Pr}\left(\left(\bigcap_{i=1}^{m} S_i \neq 0\right) \cap (A_0 \cap \cdots \cap A_{m-1})\right) \leq \mathbf{Pr}\left(\bigcap_{i=1}^{m} \tilde{S}_i \neq 0\right)$$
$$\leq \mathbf{Pr}(S_0 > \eta m) + 2e^{-\eta m},$$

as required. $\square$

We now introduce a natural coupling of $n$-queue queue-lengths processes $(X_t)$ with different initial states.

Arrival times form a Poisson process at rate $\lambda n$, and there is a corresponding sequence of uniform choices of lists of $d$ queues. Departure times form a Poisson process at rate $n$, and there is a corresponding sequence of uniform selections of a queue, except that departures from empty queues are ignored. These four processes are independent. Denote the arrival time process by $\mathbf{T}$, the choices process by $\mathbf{D}$, the departure time process by $\tilde{\mathbf{T}}$ and the selection process by $\tilde{\mathbf{D}}$.

Suppose that we are given a sequence of arrival times $\mathbf{t}$ with corresponding queue choices $\mathbf{d}$, and a sequence of departure times $\tilde{\mathbf{t}}$ with corresponding selections $\tilde{\mathbf{d}}$ of a queue, where all these times are distinct. For each possible initial queue-lengths vector $x \in \Omega = (\mathbb{Z}^+)^n$ this yields a deterministic queue-lengths process $(x_t)$ with $x_0 = x$: let us write $x_t = s_t(x; \mathbf{t}, \mathbf{d}, \tilde{\mathbf{t}}, \tilde{\mathbf{d}})$. Then for each $x \in \Omega$, the process $(s_t(x; \mathbf{T}, \mathbf{D}, \tilde{\mathbf{T}}, \tilde{\mathbf{D}}))$ has the distribution of a queue-lengths process with initial state $x$.

LEMMA 2.3. *Fix any 4-tuple $\mathbf{t}, \mathbf{d}, \tilde{\mathbf{t}}, \tilde{\mathbf{d}}$ as above, and for each $x \in \Omega$ write $s_t(x)$ for $s_t(x; \mathbf{t}, \mathbf{d}, \tilde{\mathbf{t}}, \tilde{\mathbf{d}})$. Then for each $x, y \in \Omega$, both $\|s_t(x) - s_t(y)\|_1$ and $\|s_t(x) - s_t(y)\|_\infty$ are nonincreasing; and further, if $0 \leq t < t'$ and $s_t(x) \leq s_t(y)$, then $s_{t'}(x) \leq s_{t'}(y)$.*

(We shall not need the result about $\|s_t(x) - s_t(y)\|_\infty$ in this paper, but it is convenient to record the result for use elsewhere, in particular in [8].)

PROOF OF LEMMA 2.3. Let $t_0$ be a jump time; let $x_{t_0-} = x$ and $y_{t_0-} = y$; and let $x_{t_0} = x'$ and $y_{t_0} = y'$.

Suppose that $t_0$ is an arrival time. We want to show that

$$\|x' - y'\|_1 \leq \|x - y\|_1 \tag{8}$$

and

$$\|x' - y'\|_\infty \leq \|x - y\|_\infty. \tag{9}$$

If the customer joins the same queue in the two processes, then of course $x' - y' = x - y$, and hence (8) and (9) hold. Suppose that the customer joins queue $i$ in the $x$-process and joins queue $j$ in the $y$-process, where $i \neq j$. Then $\delta_x = x(j) - x(i) \geq 0$ and $\delta_y = y(i) - y(j) \geq 0$; and by the tie-breaking rule $\delta_x + \delta_y > 0$.

Suppose first that (8) does not hold. Then we must have $x(i) \geq y(i)$ and $y(j) \geq x(j)$, and so

$$x(i) \geq y(i) = y(j) + \delta_y \geq x(j) + \delta_y = x(i) + \delta_x + \delta_y > x(i),$$

a contradiction. Hence (8) must hold.

Now suppose that (9) does not hold. Then either $x(i) - y(i) = \|x - y\|_\infty$ or $y(j) - x(j) = \|x - y\|_\infty$. But we cannot have $x(i) - y(i) = \|x - y\|_\infty$, since

$$x(j) - y(j) = x(i) + \delta_x - (y(i) - \delta_y) > x(i) - y(i).$$

Similarly we cannot have $y(j) - x(j) = \|x - y\|_\infty$, and so (9) must hold.

Suppose now that $t_0$ is a departure time, from queue $i$. If both queues are nonempty or both are empty, then $x' - y' = x - y$, and hence of course (8) and (9) hold. If exactly one queue is nonempty, then $|x'_i - y'_i| = |x_i - y_i| - 1$, and so again (8) and (9) hold.

The final comment on monotonicity is straightforward. For consider a jump time $t_0$ with $x, x', y$ and $y'$ defined as above, and suppose that $x \leq y$. If $t_0$ is a departure time, then clearly $x' \leq y'$, so suppose that $t_0$ is an arrival time. But if the new customer joins queue $i$ in the $x$-process and if $x(i) = y(i)$, then the customer joins queue $i$ also in the $y$-process, so $x' \leq y'$. □

The *position* of a customer refers to first-in–first-out queue discipline; that is, for a given customer in queue $j$ at time $t$, her position is one plus the number of customers in queue $j$ at time $t$ who arrived before her. Given a queue-lengths vector $x$ and a nonnegative integer $i$, let $\ell(i, x)$ be the number of queues with length at least $i$. We shall be interested in $\ell(i, X_t)$, the random number of queues with length at least $i$ at time $t$. Observe that if $\|x\|_1 \leq cn$, then $\ell(i, x) \leq cn/i$: this is how we shall ensure that $\ell(i, X_t)$ is not too large.

Next, let us consider the equilibrium distribution, and note some upper bounds on the total number of customers in the system and on the maximum queue length, which follow from the easy case $d = 1$.

LEMMA 2.4. (a) *For any constant $c > \frac{\lambda}{1-\lambda}$, there is a constant $\eta > 0$ such that for each positive integer $n$, in equilibrium the queue-lengths process $(X_t)$ satisfies*

$$\mathbf{Pr}(\|X_t\|_1 > cn) \leq e^{-\eta n}$$

*for each time $t \geq 0$.*

(b) *For each positive integer $n$, in equilibrium the maximum queue length $M_t$ satisfies*

$$\mathbf{Pr}(M_t \geq k) \leq n\lambda^k$$

*for each positive integer $k$ and each time $t \geq 0$.*

PROOF. For both parts of the lemma, it suffices to consider the case $d = 1$; for, as follows from a coupling result in [16] (see also [3]), if $d \leq d'$, then in equilibrium for each $k$ the number of customers with position at least $k$ with $d'$ choices is stochastically at most the corresponding number with $d$ choices. (Note that the maximum queue length $M_t$ is at least $k$ if and only if at time $t$ there is at least one customer with position at least $k$.) So, suppose that $d = 1$.

But now by the splitting property of the Poisson process the $n$ queue lengths $X_t(j)$ are independent; and each has the geometric distribution where $\mathbf{Pr}(X_t(j) = k) = (1 - \lambda)\lambda^k$, with mean $\lambda/(1 - \lambda)$. Thus the total number of customers is a sum of $n$ i.i.d. random variables with finite moment-generating function in some neighborhood of 0, and part (a) follows easily. For part (b), note that

$$\mathbf{Pr}(M_t \geq k) \leq n\mathbf{Pr}(X_t(1) \geq k) = n\lambda^k. \qquad \square$$

The bound in part (a) above extends easily over time.

LEMMA 2.5. *Let $c > \frac{\lambda}{1-\lambda}$ be a constant. Then there is a constant $\eta > 0$ such that for each positive integer $n$, in equilibrium the queue-lengths process $(X_t)$ satisfies*

$$\mathbf{Pr}(\|X_t\|_1 > cn \text{ for some } t \in [0, e^{\eta n}]) \leq 2e^{-\eta n}.$$

PROOF. Let $\frac{\lambda}{1-\lambda} < c' < c$ and let $\varepsilon = c - c' > 0$. By part (a) of the last result, there exists a constant $\eta > 0$ such that for each positive integer $n$ and each $t \geq 0$ we have $\mathbf{Pr}(\|X_t\|_1 > c'n) \leq e^{-3\eta n}$; and we may assume that $\eta < \varepsilon/18$.

Let $\delta = \varepsilon/2\lambda$. Let $j = \lceil e^{\eta n}/\delta \rceil$, and consider times $t_r = r\delta$, for $r = 0, \ldots, j$. The mean number of arrivals in a subinterval $[t_{r-1}, t_r)$ of length $\delta$ is $\varepsilon n/2$, so by (5) the probability that more than $\varepsilon n$ arrivals occur is at most $e^{-\varepsilon n/6}$. Then

$$\mathbf{Pr}(\|X_t\|_1 > cn \text{ for some } t \in [0, e^{\eta n}])$$
$$\leq \sum_{r=0}^{j} \mathbf{Pr}(\|X_{t_r}\|_1 > c'n) + \sum_{r=1}^{j} \mathbf{Pr}([t_{r-1}, t_r) \text{ has } > \varepsilon n \text{ arrivals})$$
$$\leq (e^{\eta n}/\delta + 2)(e^{-3\eta n} + e^{-\varepsilon n/6}) \leq e^{-\eta n}$$

provided $n$ is sufficiently large, and the lemma follows (for a suitable new value of $\eta$). □

We say that two states are *adjacent* if they differ by adding one customer to one of the queues in one of the states. The following lemma shows that two queue-lengths processes $(X_t)$ and $(X_t')$ will coalesce rapidly if $X_0$ and $X_0'$ are adjacent, the "unbalanced" queues are not too long and the total numbers of customers are not too large.

First we fix some constants. Let $0 < \lambda < 1$ and let $d$ be a positive integer. Let $c > \frac{\lambda}{1-\lambda}$, and let $\eta > 0$ be as in Lemma 2.5. (This will cause no loss of generality, as the case $c > \frac{\lambda}{1-\lambda}$ of course implies the case $c \leq \frac{\lambda}{1-\lambda}$ in Theorem 1.1.) Let $\varepsilon > 0$ satisfy $\nu = d\lambda\varepsilon^{d-1} < 1$. Let $k_0 = \lceil 2c/\varepsilon \rceil$. We shall keep all these constants fixed from now on, until the end of the section.

LEMMA 2.6. *There exist constants $\alpha, \beta > 0$ such that the following holds. Let $n$ be a positive integer. Let $x, x'$ be a pair of states such that $x(k) = x'(k) - 1$, $x(j) = x(j)$ for all $j \neq k$, and $\|x'\|_1 \leq cn$. Consider the queue-lengths processes $(X_t)$ and $(X_t')$, given that $X_0 = x$ and $X_0' = x'$. For all times $t \geq \alpha\|x'\|_\infty$ we have*

$$\mathbf{E}\|X_t - X_t'\|_1 = \mathbf{Pr}(X_t \neq X_t') \leq e^{-\beta t} + 2e^{-\beta n}.$$

PROOF. By Lemma 2.3, $X_t$ and $X_t'$ are always either neighbors or equal, always $X_t \leq X_t'$, and if for some time $s$ we have $X_s = X_s'$, then $X_t = X_t'$ for all $t \geq s$. Thus in particular $\mathbf{E}\|X_t - X_t'\|_1 = \mathbf{Pr}(X_t \neq X_t')$.

Initially, the queue $k$ is "unbalanced" [i.e., $X_0(k) \neq X_0'(k)$] and all other queues are "balanced." Observe that the index of the unbalanced queue in the coupled process may change over time. [E.g., suppose that $d = 2$, and just before an arrival time $t$, queue $i$ is unbalanced, and queues $i$ and $j$ are chosen. Suppose further that $X'_{t-}(i) = X'_{t-}(j)$ or $X_{t-}(i) = X_{t-}(j)$. In the former case we select queue $i$ for the $(X_t)$ process, but we will select queue $j$ for the $(X'_t)$ process if we chose $j$ before $i$. In the latter case we select queue $j$ for the $(X'_t)$ process, but we will select queue $i$ for the $(X_t)$ process if we chose $i$ before $j$. In both cases, it will now be queue $j$ that is unbalanced.]

Let $W_t$ denote the longer of the unbalanced queue lengths at time $t$ if there is such a queue, and let $W_t = 0$ otherwise. The time for the two processes to coalesce is the time $T$ until $W_t$ hits 0. We shall use Lemma 2.2 to give a suitable upper bound on $\mathbf{Pr}(W_t > 0)$. The idea is that with high probability the total number of customers in the system is not too large, hence the unbalanced queue length $W_t$ will often be driven below $k_0$, and then there is a chance of going all the way down to 0.

For each time $t \geq 0$ let $B_t$ be the event that $\|X'_s\|_1 \leq 2cn$ for each $s \in [0, t)$. It follows from Lemmas 2.3 and 2.5 that there is a constant $\eta > 0$ such that $\mathbf{Pr}(\overline{B_t}) \leq 2e^{-\eta n}$ for each positive integer $n$ and each time $t \in [0, t_0]$, where $t_0 = e^{\eta n}$. To see this, note that if $(X''_t)$ is a copy of the process such that $X''_0 = \mathbf{0}$ a.s., then there is a coupling such that $\|X'_t - X''_t\|_1 \leq cn$ for all times $t \geq 0$. But then we can couple $(X''_t)$ with an equilibrium process $(\hat{X}_t)$ so that $X''_t \leq \hat{X}_t$ for all times $t \geq 0$, and thus $\mathbf{Pr}(\|X''_t\|_1 > cn$ for some $t \in [0, t_0]) \leq 2e^{-\eta n}$. Finally, if $\|X''_t\|_1 \leq cn$ and $\|X'_t - X''_t\|_1 \leq cn$, then $\|X'_t\|_1 \leq 2cn$.

We need some notation concerning the jumps in the unbalanced queue length $W_t$. Let $N_t$ be the number of such jumps in the interval $[0, t]$. Also let $N = N_T$, the total number of these jumps. Let $T_j$ be the time of the $j$th jump if $N \geq j$, and otherwise let $T_j$ be the coalescence time $T$.

Let $S_0 = x'(k) = W_0$, the longer unbalanced queue length at time $t = 0$. For each positive integer $j$, if $N \geq j$, let $S_j = W_{T_j}$, which is either 0 or the longer of the unbalanced queue lengths at time $T_j$, immediately after the $j$th arrival or departure at the unbalanced queue. Also, if $N \geq j$, let $Z_j$ be the $\pm 1$-valued random variable $S_j - S_{j-1}$. For each nonnegative integer $j$, let $\phi_j$ be the $\sigma$-field generated by all events before time $T_{j+1}$. Let also $A_j$ be the $\phi_j$-measurable event $B_{T_{j+1}}$.

Let
$$\delta = \min\left\{\frac{1}{d\lambda + 1}, \frac{1}{\nu + 1} - \frac{1}{2}\right\},$$
and note that $\delta > 0$ since $0 < \nu < 1$. We shall use Lemma 2.2, with this value of $\delta$. Note first that the arrival rate at the longer of the unbalanced queues is always at most $d\lambda$, and the departure rate is 1. Thus on the event $N \geq j$ we have
$$\mathbf{Pr}(Z_j = -1|\phi_{j-1}) \geq \frac{1}{d\lambda + 1} \geq \delta.$$

The key observation is that, on the event $\{N \geq j\} \cap A_{j-1} \cap \{S_{j-1} \geq k_0\}$, we have
$$\mathbf{Pr}(Z_j = -1|\phi_{j-1}) \geq \frac{1}{\nu + 1} \geq \frac{1}{2} + \delta.$$

For consider a time $t > 0$. Note that on $B_t$ we have $\ell(k_0, X'_{t-}) \leq 2cn/k_0 \leq \varepsilon n$, that is, there are at most $\varepsilon n$ queues with length at least $k_0$. Suppose that $T > t$, $B_t$ holds and $W_t \geq k_0$. Arrivals into the system occur at rate $n\lambda$, and the probability that such an arrival joins the longer unbalanced queue is at most $\frac{d}{n}\varepsilon^{d-1}$. Hence the rate of arrivals at the longer unbalanced queue is at most $\nu = d\lambda\varepsilon^{d-1}$, whereas the rate of departures is 1, and so the next jump is a departure with probability at least $\frac{1}{\nu+1}$.

We have now shown that on the event $N \geq m$, $S_m - S_0$ can be written as a sum $\sum_{i=1}^m Z_i$ for $\{-1, 1\}$-valued random variables $Z_i$ that satisfy the conditions of Lemma 2.2, with the same notation for $\delta$, etc. Hence there exists a constant $\eta_1 > 0$, such that for all $m \geq m_0 = \lceil \eta_1^{-1} x'(k) \rceil$

$$\mathbf{Pr}\left(\{N \geq m\} \cap \left(\bigcap_{i=0}^{m-1} A_i\right) \cap \left(\bigcap_{i=1}^{m} \{S_i \neq 0\}\right)\right) \leq 2e^{-\eta_1 m}.$$

Let $n$ be sufficiently large that $2m_0 \leq t_0 = e^{\eta n}$, let $t$ satisfy $2m_0 \leq t \leq t_0$ and let $m = \lfloor t/2 \rfloor$. Then, since jumps occur at rate at least 1 while the queue is nonempty, by (4)

$$\mathbf{Pr}(\{T > t\} \cap \{N_t < m\}) \leq e^{-t/8}.$$

Also,

$$\mathbf{Pr}\left(\{N_t \geq m\} \cap \left(\bigcup_{i=0}^{m-1} \overline{A_i}\right)\right) \leq \mathbf{Pr}(\overline{B_t}) \leq 2e^{-\eta n}.$$

This gives the desired upper bound on $\mathbf{Pr}(T > t)$, since

$$\mathbf{Pr}(T > t) \leq \mathbf{Pr}(\{T > t\} \cap \{N_t < m\}) + \mathbf{Pr}\left(\{N_t \geq m\} \cap \left(\bigcup_{i=0}^{m-1} \overline{A_i}\right)\right)$$
$$+ \mathbf{Pr}\left(\{N_t \geq m\} \cap \left(\bigcap_{i=0}^{m-1} A_i\right) \cap \left(\bigcap_{i=1}^{m} \{S_i \neq 0\}\right)\right).$$

Thus we have shown that
$$\mathbf{Pr}(X_t \neq X'_t) \leq e^{-\beta t} + 2e^{-\beta n}$$
whenever $2m_0 \leq t \leq t_0$.

Now we check that we can drop the upper bound $t_0$ on $t$. For let $n$ be sufficiently large that $e^{-\beta t_0} \leq e^{-\beta n}$. If $t > t_0$, then
$$\mathbf{Pr}(X_t \neq X'_t) \leq \mathbf{Pr}(X_{t_0} \neq X'_{t_0}) \leq 3e^{-\beta n},$$

and so

$$\mathbf{Pr}(X_t \neq X'_t) \leq e^{-\beta t} + 3e^{-\beta n}$$

for all $t \geq 2m_0$.

Finally, by replacing $\beta$ with a smaller constant $\beta > 0$, we can replace 3 by 2, and ensure that the inequality holds for all positive integers $n$, not just for sufficiently large $n$. □

Recall that we have fixed a constant $c > \lambda/(1-\lambda)$ and that we are using the coupling introduced before Lemma 2.3.

LEMMA 2.7. *Let $\alpha$ and $\beta$ be as in the last lemma. Let $(X_t)$ and $(Y_t)$ have any initial distributions. Let $t \geq 0$, and let the "bad" set $B$ be the set of $z \in \Omega = (\mathbb{Z}^+)^n$ such that $\|z\|_1 > cn$ or $\|z\|_\infty > t\alpha^{-1}$. Then*

$$\mathbf{E}[\|X_t - Y_t\|_1 \mathbf{1}_{\{X_0 \in \overline{B}\} \cap \{Y_0 \in \overline{B}\}}] \leq 2cn(e^{-\beta t} + 2e^{-\beta n}).$$

PROOF. Given two distinct states $x$ and $y$ in $\overline{B} = \Omega \setminus B$, we can choose a path $x = z_0, z_1, \ldots, z_m = y$ of adjacent states in $\overline{B}$ from $x$ down to the all zero state and back up to $y$, where $m \leq \|x\|_1 + \|y\|_1$. Let us write $(X_t^x)$ to denote the queue-lengths process starting at $x$, and similarly for $(Y_t^y)$ (so in fact $X_t^x = Y_t^x$ always). By the last lemma, for all states $x$ and $y$ in $\overline{B}$

$$\mathbf{E}[\|X_t^x - Y_t^y\|_1] \leq \sum_{i=0}^{m-1} \mathbf{E}[\|X_t^{z_i} - X_t^{z_{i+1}}\|_1] \leq (\|x\|_1 + \|y\|_1)(e^{-\beta t} + 2e^{-\beta n}).$$

Hence the lemma follows. □

We may now complete the proofs of Theorem 1.1 and Lemma 2.1, by taking $(Y_t)$ to be in equilibrium, with $X_0$ and $Y_0$ independent, and handling the "bad" initial states appropriately. Consider first Theorem 1.1. Note that

$$d_{\mathrm{TV}}(\mathcal{L}(X_t), \mathbf{\Pi}) \leq \mathbf{Pr}(X_t \neq Y_t)$$
$$\leq \mathbf{E}[\mathbf{1}_{X_t \neq Y_t} \mathbf{1}_{\{X_0 \in \bar{B}\} \cap \{Y_0 \in \bar{B}\}}] + \mathbf{Pr}(X_0 \in B) + \mathbf{Pr}(Y_0 \in B)$$
$$\leq \mathbf{E}[\|X_t - Y_t\|_1 \mathbf{1}_{\{X_0 \in \bar{B}\} \cap \{Y_0 \in \bar{B}\}}] + \mathbf{Pr}(X_0 \in B) + \mathbf{Pr}(Y_0 \in B).$$

But $\mathbf{Pr}(\|Y_0\|_1 > cn) = e^{-\Omega(n)}$ (see Lemma 2.5), since $c > \lambda/(1-\lambda)$. Also, from (25) in Section 5, for each $j \in \{1, \ldots, n\}$ and each nonnegative integer $i$, $\mathbf{Pr}(Y_0(j) \geq i) \leq \lambda^i$; and it follows that

$$\mathbf{Pr}(\|Y_0\|_\infty > t\alpha^{-1}) \leq ne^{-\Omega(t)}.$$

Theorem 1.1 now follows from Lemma 2.7.

Finally let us complete the proof of Lemma 2.1. Note that $d_{\mathrm{W}}(\mathcal{L}(X_t), \mathbf{\Pi}) \leq \mathbf{E}[\|X_t - Y_t\|_1]$. We break $\mathbf{E}[\|X_t - Y_t\|_1]$ into the sum of two parts

$$\mathbf{E}[\|X_t - Y_t\|_1 \mathbf{1}_{\{X_0 \in \overline{B}\} \cap \{Y_0 \in \overline{B}\}}] + \mathbf{E}[\|X_t - Y_t\|_1 \mathbf{1}_{\{X_0 \in B\} \cup \{Y_0 \in B\}}].$$

The first part is bounded in Lemma 2.7. For the second, we have by Lemma 2.3 that

$$\begin{aligned}
\mathbf{E}[\|X_t &- Y_t\|_1 \mathbf{1}_{\{X_0 \in B\} \cup \{Y_0 \in B\}}] \\
&\leq \mathbf{E}[\|X_0 - Y_0\|_1 \mathbf{1}_{\{X_0 \in B\} \cup \{Y_0 \in B\}}] \\
&\leq \mathbf{E}[(\|X_0\|_1 + \|Y_0\|_1)(\mathbf{1}_{\|X_0\|_1 > cn} + \mathbf{1}_{\|X_0\|_1 \leq cn, \|Y_0\|_1 > cn})] \\
&\quad + \mathbf{E}[(\|X_0\|_1 + \|Y_0\|_1) \mathbf{1}_{\|X_0\|_1 \leq cn, \|Y_0\|_1 \leq cn, \max\{\|X_0\|_\infty, \|Y_0\|_\infty\} > t\alpha^{-1}}] \\
&\leq \mathbf{E}[\|X_0\|_1 \mathbf{1}_{\|X_0\|_1 > cn}] + \mathbf{E}[\|Y_0\|_1]\mathbf{Pr}(\|X_0\|_1 > cn) + cn\mathbf{Pr}(\|Y_0\|_1 > cn) \\
&\quad + \mathbf{E}[\|Y_0\|_1 \mathbf{1}_{\|Y_0\|_1 > cn}] + 2cn(\mathbf{Pr}(\|X_0\|_\infty > t\alpha^{-1}) + \mathbf{Pr}(\|Y_0\|_\infty > t\alpha^{-1})),
\end{aligned}$$

where the last inequality uses the independence of $X_0$ and $Y_0$. Now we may use the estimates above concerning $Y_0$, together with the fact that $\mathbf{E}[\|Y_0\|_1] \leq \frac{\lambda}{1-\lambda} n$, to complete the proof. (Note that

$$\mathbf{E}[\|Y_0\|_1]\mathbf{Pr}(\|X_0\|_1 > cn) \leq \frac{\lambda}{1-\lambda} n \mathbf{Pr}(\|X_0\|_1 > cn) \leq \mathbf{E}[\|X_0\|_1 \mathbf{1}_{\|X_0\|_1 > cn}]$$

since $c \geq \frac{\lambda}{1-\lambda}$.)

**3. One choice: proof of Theorem 1.2.** Let $0 < \lambda < 1$ and let $d = 1$. Let $X_0$ be in equilibrium.

*Part* (a). The queue lengths $X_t(j)$ for $j = 1, \ldots, n$ are independent geometric random variables with parameter $\lambda$, and so $\mathbf{Pr}(X_t(j) \geq m) = \lambda^m$ for each nonnegative integer $m$. Hence

$$\mathbf{Pr}(M_t \leq m) = (1 - \lambda^{m+1})^n,$$

and so

(10) $\quad \exp(-n\lambda^{m+1}/(1 - \lambda^{m+1})) \leq \mathbf{Pr}(M_t \leq m) \leq \exp(-n\lambda^{m+1}).$

The rest of part (a) follows easily.

*Part* (b). The proofs we give for parts (b) and (c) are similar to the proofs of the corresponding parts of Theorem 1.2 in [7]. Let $V = \min_{0 \leq t \leq \tau} M_t$. Let $\varepsilon > 0$, let $m = m(n) = \lfloor (1-\varepsilon) \frac{\ln n}{\ln(1/\lambda)} \rfloor$ and let $\theta = \frac{1}{n^2\tau}$. Consider the times $i\theta$ for $i = 0, 1, \ldots, \lceil \tau/\theta \rceil$. Let $A$ be the event that $M_{i\theta} < m$ for some $i \in \{0, 1, \ldots, \lceil \tau/\theta \rceil\}$. We have $n\lambda^m = \Omega(n^\varepsilon)$, and hence by (10)

$$\mathbf{Pr}(A) \leq (\tau/\theta + 1)\exp(-n\lambda^m) = o(1).$$

Let $B$ be the event that some queue receives at least two customers in some time interval $[i\theta, (i+1)\theta)$, where $i \in \{0, 1, \ldots, \lceil \tau/\theta \rceil\}$. For each queue, the number of arrivals in the interval $[i\theta, (i+1)\theta)$ is $Po(\lambda\theta)$, and so the probability that there are at least two arrivals is at most $(\lambda\theta)^2$. Hence

$$\mathbf{Pr}(B) \leq (\tau/\theta + 1)n(\lambda\theta)^2 = O(n\tau\theta) = o(1).$$

But

$$\mathbf{Pr}(V \leq m - 2) \leq \mathbf{Pr}(A) + \mathbf{Pr}(B),$$

and so $V \geq m - 1$ a.a.s. But by part (a)

$$V \leq M_0 \leq (1+\varepsilon)\frac{\ln n}{\ln(1/\lambda)} \quad \text{a.a.s.,}$$

and part (b) follows.

*Part* (c). Let $Z = \max_{0 \leq t \leq n^K} M_t$. Let $\varepsilon > 0$. We show first that

(11) $$Z \leq (K + 1 + \varepsilon)\frac{\ln n}{\ln(1/\lambda)} \quad \text{a.a.s.}$$

We argue much as in the proof of part (b). Let $\theta = \exp(-\ln n/\ln \ln n)$. Consider the times $i\theta$ for $i = 0, 1, \ldots, \lceil n^K/\theta \rceil$. Let $k = k(n) = \lceil (K + 1 + \varepsilon/2) \ln n/\ln(1/\lambda) \rceil$, and let $A$ be the event that $M_{i\theta} \geq k$ for some $i \in \{0, 1, \ldots, \lceil n^K/\theta \rceil\}$. Then since $\mathbf{Pr}(M_0 \geq k) \leq n\lambda^k$,

$$\begin{aligned}
\mathbf{Pr}(A) &\leq (n^K/\theta + 1)n\lambda^k \\
&= \exp((K + 1 + o(1))\ln n - (K + 1 + \varepsilon/2 + o(1))\ln n) \\
&= \exp(-(\varepsilon/2 + o(1))\ln n) \\
&\to 0 \quad \text{as } n \to \infty.
\end{aligned}$$

For each queue, the number of arrivals in the interval $[i\theta, (i+1)\theta)$ is $Po(\lambda\theta)$, and

$$\begin{aligned}
\mathbf{Pr}(Po(\lambda\theta) \geq \ln n/(\ln \ln n)^2) &\leq (\lambda\theta)^{\ln n/(\ln \ln n)^2} \\
&= \exp(-\Omega(\ln^2 n/(\ln \ln n)^3)).
\end{aligned}$$

Thus, if $B$ is the event that some queue in some interval $[i\theta, (i+1)\theta)$ where $i \in \{0, 1, \ldots, \lceil n^K/\theta \rceil\}$ receives at least $\ln n/(\ln \ln n)^2$ customers, then

$$\mathbf{Pr}(B) = \exp(-\Omega(\ln^2 n/(\ln \ln n)^3)).$$

But, for $n$ sufficiently large,

$$\mathbf{Pr}\left(Z > (K+1+\varepsilon)\frac{\ln n}{\ln(1/\lambda)}\right) \leq \mathbf{Pr}(A) + \mathbf{Pr}(B) \to 0 \quad \text{as } n \to \infty,$$

and so (11) holds.

Now let $0 < \varepsilon < 1$, and let $k = k(n) = \lceil (K + 1 - \varepsilon) \ln n / \ln(1/\lambda) \rceil$. We will show that $Z \geq k$ a.a.s., which will complete the proof of this part and thus of the theorem. For each time $t > 0$, let $\phi_t$ be the $\sigma$-field generated by all events until time $t$. Let $c > \lambda/(1-\lambda)$. Let $C$ be the event that $\|X_t\|_1 \leq cn$ for each $t \in [0, n^K]$. Then $C$ holds a.a.s. by Lemma 2.5. Also, by Theorem 1.1 there are constants $n_0$ and $\eta > 0$ such that the following holds. Let $n \geq n_0$ and consider the system with $n$ queues. Let $x$ be a queue-lengths vector such that $\|x\|_1 \leq cn$ and $\|x\|_\infty \leq k - 1$. Then, given $X_0 = x$,
$$d_{\mathrm{TV}}(\mathcal{L}(X_t), \Pi) \leq e^{-\ln^2 n}$$
for all times $t \geq t_1 = \eta^{-1} \ln^2 n$. In particular, for $n$ sufficiently large by (10)
$$\mathbf{Pr}(M_{t_1} \leq k - 1 | X_0 = x) \leq e^{-n\lambda^k} + e^{-\ln^2 n}.$$
Thus, since the system is in equilibrium, for each $i = 0, 1, \ldots$
$$\mathbf{Pr}(M_{(i+1)t_1} \leq k - 1 | \phi_{it_1}) \leq e^{-n\lambda^k} + e^{-\ln^2 n}$$
on the event $D_i = \{\|X_{it_1}\|_1 \leq cn\} \cap \{M_{it_1} \leq k - 1\}$. Hence if we denote $\lfloor n^K / t_1 \rfloor$ by $i_0$,

$$\begin{aligned}
\mathbf{Pr}(\{Z \leq k-1\} \cap C) &\leq \mathbf{Pr}\left(\bigcap_{i=0}^{i_0} D_i\right) \\
&\leq \mathbf{Pr}(D_0) \prod_{i=0}^{i_0-1} \mathbf{Pr}\left(D_{i+1} \Big| \bigcap_{j=0}^{i} D_j\right) \\
&\leq (e^{-n\lambda^k} + e^{-\ln^2 n})^{i_0} \\
&\leq \exp(-(n^K/t_1) n^{-(K-\varepsilon+o(1))}) \\
&= \exp(-n^{\varepsilon+o(1)}) \\
&\to 0 \quad \text{as } n \to \infty.
\end{aligned}$$

Since $C$ holds a.a.s., it follows that $Z \geq k$ a.a.s., as required.

**4. Concentration.** In this section we prove concentration of measure results for the queue-lengths process $(X_t)$. Let $n$ be a positive integer, and let $\Omega$ be the corresponding set of queue-lengths vectors, that is, the set of nonnegative vectors in $\mathbb{Z}^n$. Let us say that a real-valued function $f$ on $\Omega$ is *Lipschitz* (with constant 1) if
$$|f(x) - f(y)| \leq \|x - y\|_1$$
for all $x, y \in \Omega$. Let $d$ be a fixed positive integer. The key result is the following lemma.

LEMMA 4.1. *There is a constant $c > 0$ such that the following holds. Let $n \geq 2$ be an integer and consider the $n$-queue system. Let the queue-lengths vector $Y$ have the equilibrium distribution. Let $f$ be a Lipschitz function on $\Omega$. Then for each $u \geq 0$*

$$\mathbf{Pr}(|f(Y) - \mathbf{E}[f(Y)]| \geq u) \leq n e^{-cu/n^{1/2}}.$$

Recall that $\ell(k, x)$ denotes $|\{j : x(j) \geq k\}|$, the number of queues of length at least $k$; and observe that for any fixed $k$ this is Lipschitz as a function of $x$. We deduce from the last lemma the following result concerning the random variables $\ell(k, Y)$.

LEMMA 4.2. *Consider the $n$-queue system, and let the queue-lengths vector $Y$ have the equilibrium distribution. For each nonnegative integer $k$ let $\ell(k) = \mathbf{E}[\ell(k, Y)]$. Then for any constant $c > 0$,*

$$\mathbf{Pr}\left(\sup_k |\ell(k, Y) - \ell(k)| \geq c n^{1/2} \ln^2 n\right) = e^{-\Omega(\ln^2 n)}.$$

*Also, for each integer $r \geq 2$*

$$\sup_k |\mathbf{E}[\ell(k, Y)^r] - \ell(k)^r| = O(n^{r-1} \ln^2 n).$$

PROOF. We argue as in the proof of Lemma 5.2 of [7]. For the first part, let $c_1 > \frac{\lambda}{1-\lambda}$, and note that by Lemma 2.5

$$\mathbf{Pr}(\ell(\lceil c_1 n \rceil, Y) > 0) = e^{-\Omega(n)}.$$

It follows that we may restrict attention to queue lengths $k < c_1 n$. For, since always $\ell(k, Y) \leq n$ we have $\ell(\lceil c_1 n \rceil) < 1$ for $n$ sufficiently large; and then

$$\mathbf{Pr}\left(\sup_{k \geq c_1 n} |\ell(k, Y) - \ell(k)| \geq 1\right) \leq \mathbf{Pr}(\ell(\lceil c_1 n \rceil, Y) \geq 1) = e^{-\Omega(n)}.$$

Now the first part of the lemma follows easily from Lemma 4.1.

For the second part, fix an integer $r \geq 2$. By Lemma 4.1 there is a constant $c_2 > 0$ such that, if we set $u = c_2 n^{1/2} \ln n$, then

$$\sup_k \mathbf{Pr}(|\ell(k, Y) - \ell(k)| > u) = o(n^{-r}).$$

Hence, for each positive integer $s \leq r$,

$$\mathbf{E}[|\ell(k, Y) - \ell(k)|^s] \leq u^s + n^s \mathbf{Pr}(|\ell(k, Y) - \ell(k)| > u) = u^s + o(1),$$

uniformly over $k$. The result now follows from

$$
\begin{aligned}
0 &\leq \mathbf{E}[\ell(k,Y)^r] - \ell(k)^r \\
&= \sum_{s=2}^{r} \mathbf{E}[(\ell(k,Y) - \ell(k))^s]\ell(k)^{r-s} \\
&\leq \sum_{s=2}^{r} \mathbf{E}[|\ell(k,Y) - \ell(k)|^s]n^{r-s} \\
&= O(n^{r-1}\ln^2 n),
\end{aligned}
$$

uniformly over $k$. $\square$

Our proof of Lemma 4.1 will follow the lines of the proof of Lemma 5.1 in [7]. The task is somewhat easier here, and we obtain tighter bounds, since for each fixed $n$ the departures process has a bounded rate. Departures occur as events in a Poisson process at rate $n$, each one from an independently selected uniformly random queue, except that departures from empty queues are ignored. As in [7], along the way we prove concentration for Lipschitz functions of the time-dependent process for "nice" initial conditions—see Lemma 4.3 below.

An overview of the proof is as follows. Consider a queue-lengths process $(X_t)$ where $X_0 = \mathbf{0}$. For $t > 0$, let $Z_t$ be the total number of arrivals in $[0,t]$, and let $\tilde{Z}_t$ be the total number of departures in $[0,t]$ (including "virtual ones," i.e., departures from empty queues). Thus $Z_t \sim Po(\lambda nt)$ and $\tilde{Z}_t \sim Po(nt)$. Let $\mu_t = \mathbf{E}[f(X_t)]$, and $\mu_t(z,\tilde{z}) = \mathbf{E}[f(X_t)|Z_t = z, \tilde{Z}_t = \tilde{z}]$. We use earlier coupling results and the bounded differences method to upper bound $\mathbf{Pr}(|f(X_t) - \mu_t(z,\tilde{z})| \geq u|Z_t = z, \tilde{Z}_t = \tilde{z})$.

Next we remove the conditioning on $Z_t$ and $\tilde{Z}_t$. To do this, we choose suitable "widths" $w$ and $\tilde{w}$, and use the fact that $\mathbf{Pr}(|Z_t - \lambda nt| > w)$ and $\mathbf{Pr}(|\tilde{Z}_t - nt| > \tilde{w})$ are small, and for $z, \tilde{z}$ such that $|z - \lambda nt| \leq w, |\tilde{z} - nt| \leq \tilde{w}$ the difference $|\mu_t(z,\tilde{z}) - \mu_t|$ is at most about $4(w + \tilde{w})$. We thus find that $\mathbf{Pr}(|f(X_t) - \mu_t| \geq 5(w + \tilde{w}))$ is small. The part of the proof up to here is encapsulated in Lemma 4.3 below. Finally we use the mixing results, Theorem 1.1 and Lemma 2.1, to relate the distribution of $X_t$ to the equilibrium distribution.

Let us start on the details of the proof. In this section we shall use the following lemma with $x_0 = \mathbf{0}$; we consider more general initial states for later use.

LEMMA 4.3. *There is a constant $c > 0$ such that the following holds. Let $n \geq 2$ be an integer and let $f$ be a Lipschitz function on $\Omega$. Let also $x_0 \in \Omega$* *and assume that the queue-lengths process* $(X_t)$ *satisfies* $X_0 = x_0$ *a.s. Then for all times* $t > 0$ *and all* $u \geq 0$,

$$\mathbf{Pr}(|f(X_t) - \mu_t| \geq u) \leq n e^{-cu^2/(nt+u)}. \tag{12}$$

PROOF. Note first that we may assume without loss of generality that $f(x_0) = 0$, and so $|f(X_t)| \leq Z_t + \tilde{Z}_t$, since we could replace $f(x)$ by its translation $\tilde{f}(x) = f(x) - f(x_0)$. Let $z$ and $\tilde{z}$ be positive integers, and condition on $Z_t = z$ and $\tilde{Z}_t = \tilde{z}$. Now $f(X_t)$ depends on $2(z + \tilde{z})$ independent random variables $T_1, \ldots, T_z, D_1, \ldots, D_z, \tilde{T}_1, \ldots, \tilde{T}_{\tilde{z}}, \tilde{D}_1, \ldots, \tilde{D}_{\tilde{z}}$ which specify the arrival time and corresponding choice of $d$ queues for each of $z$ customers, and the $\tilde{z}$ departure times and corresponding selections of a queue during $[0, t]$. This property relies on the well-known fact that conditional on the number of events of a Poisson process during $[0, t]$, the arrival times are a sample of i.i.d. random variables uniform on $[0, t]$. Let $\mathbf{T} = (T_1, \ldots, T_z)$, $\mathbf{D} = (D_1, \ldots, D_z)$, $\tilde{\mathbf{T}} = (\tilde{T}_1, \ldots, \tilde{T}_{\tilde{z}})$ and $\tilde{\mathbf{D}} = (\tilde{D}_1, \ldots, \tilde{D}_{\tilde{z}})$. We may write $f(X_t)$ as $g(\mathbf{T}, \mathbf{D}, \tilde{\mathbf{T}}, \tilde{\mathbf{D}})$ where

$$g(\mathbf{t}, \mathbf{d}, \tilde{\mathbf{t}}, \tilde{\mathbf{d}}) = f(s_t(x_0, \mathbf{t}, \mathbf{d}, \tilde{\mathbf{t}}, \tilde{\mathbf{d}})),$$

in the notation of Lemma 2.3.

We now prove that, conditional on $Z_t = z$ and $\tilde{Z}_t = \tilde{z}$, the random variable $f(X_t)$ is strongly concentrated, by showing that $g(\mathbf{t}, \mathbf{d}, \tilde{\mathbf{t}}, \tilde{\mathbf{d}})$ satisfies a "bounded differences" condition. Suppose first that we alter a single coordinate value $d_j$ or $\tilde{d}_j$. Then the value of $g$ can change by at most 2; by Lemma 2.3 starting at time $t_j$ with $\|x_{t_j} - y_{t_j}\|_1 \leq 2$. Similarly, if we change a coordinate value $t_j$ or $\tilde{t}_j$, the value of $g$ can change by at most 2; we may see this by applying Lemma 2.3 once at the earlier time and once at the later time. Now we use the independent bounded differences inequality; see, for instance, [12]. Hence, for each $u > 0$

$$\mathbf{Pr}(|g(\mathbf{T}, \mathbf{D}, \tilde{\mathbf{T}}, \tilde{\mathbf{D}}) - \mathbf{E}[g(\mathbf{T}, \mathbf{D}, \tilde{\mathbf{T}}, \tilde{\mathbf{D}})]| \geq u) \leq 2\exp\left(-\frac{u^2}{4(z+\tilde{z})}\right).$$

In other words, we have proved that for any $u > 0$

$$\mathbf{Pr}(|f(X_t) - \mu_t(z, \tilde{z})| \geq u | Z_t = z, \tilde{Z}_t = \tilde{z}) \leq 2\exp\left(-\frac{u^2}{4(z+\tilde{z})}\right), \tag{13}$$

which is the desired upper bound on the quantity on the left-hand side.

Next we will remove the conditioning on $Z_t$. We will choose suitable "widths" $w = w(n) > 0$ and $\tilde{w} = \tilde{w}(n) > 0$, where $0 \leq w \leq \lambda nt$ and $0 \leq \tilde{w} \leq nt$. Let $I$ denote the interval of integer values $z$ such that $|z - \lambda nt| \leq w$; let $\tilde{I}$ denote the interval of integer values $\tilde{z}$ such that $|\tilde{z} - nt| \leq \tilde{w}$. Since $Z_t \sim Po(\lambda nt)$ and $\tilde{Z}_t \sim Po(nt)$, by inequalities (4) and (5),

$$\mathbf{Pr}(Z_t \notin I) = \mathbf{Pr}(|Z_t - \lambda nt| > w) \leq 2\exp\left(-\frac{w^2}{3\lambda nt}\right) \tag{14}$$

and

$$\textbf{Pr}(\tilde{Z}_t \notin \tilde{I}) = \textbf{Pr}(|\tilde{Z}_t - nt| > \tilde{w}) \leq 2\exp\left(-\frac{\tilde{w}^2}{3nt}\right). \tag{15}$$

We shall assume that $w \geq 2(\lambda nt \ln n)^{1/2}$ and $\tilde{w} \geq 2(nt \ln n)^{1/2}$, and so it follows by (14) and (15) that

$$\mathbf{E}[(Z_t + \tilde{Z}_t)\mathbf{1}_{\{Z_t \notin I\} \cup \{\hat{Z}_t \notin \tilde{I}\}}] = o(1). \tag{16}$$

From Lemma 2.3, for each $z$

$$|\mu_t(z+1, \tilde{z}) - \mu_t(z, \tilde{z})| \leq 1 \tag{17}$$

and

$$|\mu_t(z, \tilde{z}+1) - \mu_t(z, \tilde{z})| \leq 1. \tag{18}$$

We claim that for each $z \in I$ and $\tilde{z} \in \tilde{I}$,

$$|\mu_t(z, \tilde{z}) - \mu_t| \leq 2(w + \tilde{w}) + o(1). \tag{19}$$

To prove this, observe that

$$\mu_t = \sum_{z \in I, \tilde{z} \in \tilde{I}} \mu_t(z, \tilde{z}) \textbf{Pr}(Z_t = z, \tilde{Z}_t = \tilde{z}) + \mathbf{E}[f(X_t)\mathbf{1}_{\{Z_t \notin I\} \cup \{\tilde{Z}_t \notin \tilde{I}\}}].$$

Hence by (16), since $|f(X_t)| \leq Z_t + \tilde{Z}_t$,

$$\mu_t \leq \max_{z \in I, \tilde{z} \in \tilde{I}}\{\mu_t(z, \tilde{z})\} + \mathbf{E}[(Z_t + \tilde{Z}_t)\mathbf{1}_{\{Z_t \notin I\} \cup \{\tilde{Z}_t \notin \tilde{I}\}}]$$
$$\leq \max_{z \in I, \tilde{z} \in \tilde{I}}\{\mu_t(z, \tilde{z})\} + o(1),$$

and by (14) and (15),

$$\mu_t \geq \min_{z \in I, \tilde{z} \in \tilde{I}}\{\mu_t(z, \tilde{z})\}\textbf{Pr}(Z_t \in I, \tilde{Z}_t \in \tilde{I}) + o(1)$$
$$\geq \min_{z \in I, \tilde{z} \in \tilde{I}}\{\mu_t(z, \tilde{z})\} + o(1).$$

Now we may use (17) and (18) to complete the proof of (19). By (13), (14), (15) and (19),

$$\mathbf{Pr}(|f(X_t) - \mu_t| \geq 3(w + \tilde{w}))$$
$$\leq \sum_{z \in I, \tilde{z} \in \tilde{I}} \mathbf{Pr}(|f(X_t) - \mu_t| \geq 3(w + \tilde{w})|Z_t = z, \tilde{Z}_t = \tilde{z})\mathbf{Pr}(Z_t = z, \tilde{Z}_t = \tilde{z})$$
$$+ \mathbf{Pr}(Z_t \notin I) + \mathbf{Pr}(\tilde{Z}_t \notin \tilde{I})$$
$$\leq \sum_{z \in I, \tilde{z} \in \tilde{I}} \mathbf{Pr}(|f(X_t) - \mu_t(z, \tilde{z})| \geq w + \tilde{w} + o(1)|Z_t = z, \tilde{Z}_t = \tilde{z})$$
$$\times \mathbf{Pr}(Z_t = z, \tilde{Z}_t = \tilde{z})$$
$$+ \mathbf{Pr}(|Z_t - \lambda nt| > w) + \mathbf{Pr}(|\tilde{Z}_t - nt| > \tilde{w})$$
$$\leq \exp\left(-\frac{(w + \tilde{w} + o(1))^2}{4(\lambda nt + nt + w + \tilde{w})}\right) + 2\exp\left(-\frac{w^2}{3\lambda nt}\right) + 2\exp\left(-\frac{\tilde{w}^2}{3nt}\right)$$
$$\leq \exp\left(-\frac{(1 + o(1))(w + \tilde{w})^2}{8(\lambda + 1)nt}\right) + 2\exp\left(-\frac{w^2}{3\lambda nt}\right) + 2\exp\left(-\frac{\tilde{w}^2}{3nt}\right).$$

Now let $u$ and $t$ satisfy

$$(20) \qquad 12(nt \ln n)^{1/2} \leq u \leq 6\lambda nt.$$

Let $w = \tilde{w} = u/6$. Then $u = 3(w + \tilde{w})$; and $w$ and $\tilde{w}$ are as required, that is, $2(\lambda nt \ln n)^{1/2} \leq w \leq \lambda nt$ and $2(nt \ln n)^{1/2} \leq \tilde{w} \leq nt$. Hence for $n$ sufficiently large we have

$$(21) \qquad \mathbf{Pr}(|f(X_t) - \mu_t| \geq u) \leq \exp\left(-\frac{u^2}{144nt}\right).$$

But if $u \leq 12(nt \ln n)^{1/2}$, then $\exp(-\frac{u^2}{144nt}) \geq n^{-1}$. Thus, as long as $u \leq 6\lambda nt$ we have

$$\mathbf{Pr}(|f(X_t) - \mu_t| \geq u) \leq n \exp\left(-\frac{u^2}{144nt}\right).$$

Now let us get rid of the upper bound on $u$. If $6\lambda nt < u < 6ent$, then by the above of course

$$\mathbf{Pr}(|f(X_t) - \mu_t| \geq u) \leq n \exp\left(-\frac{\lambda^2}{144e^2}\frac{u^2}{nt}\right).$$

Finally consider $u \geq 6ent$. We saw that $|f(X_t)| \leq Z_t + \tilde{Z}_t$. Thus $|\mu_t| \leq \mathbf{E}[|f(X_t)|] \leq (1 + \lambda)nt \leq 2nt$. Hence, if $u \geq 6ent$, then

$$\mathbf{Pr}(|f(X_t) - \mu_t| \geq u) \leq \mathbf{Pr}(|f(X_t)| \geq 2u/3)$$
$$\leq \mathbf{Pr}(Z_t \geq u/3) + \mathbf{Pr}(\tilde{Z}_t \geq u/3)$$
$$\leq 2\mathbf{Pr}(\tilde{Z}_t \geq u/3).$$

But $\tilde{Z}_t \sim Po(nt)$ and $u/3 \geq 2ent$, and so by (6) the last bound is at most $2^{1-u/3}$. The lemma now follows. $\square$

We shall use Lemma 4.3 here with $X_0 = \mathbf{0}$ to complete the proof of Lemma 4.1. As we saw before, we may assume that $f(\mathbf{0}) = 0$, and hence always $|f(x)| \leq \|x\|_1$. It remains to relate the distribution of $X_t$ with $X_0 = \mathbf{0}$ to the equilibrium distribution. But by Theorem 1.1 there exists a constant $\eta > 0$ such that, for each positive integer $n$ and each time $t \geq 0$, if $Y$ has the equilibrium distribution, then we have $d_{\mathrm{TV}}(\mathcal{L}(X_t), \mathcal{L}(Y)) \leq ne^{-\eta t} + 2e^{-\eta n}$. Also, by Lemma 4.3 we may assume that $\eta > 0$ is sufficiently small that, for each $n$, each $t > 0$ and each $u \geq 0$

$$\mathbf{Pr}(|f(Y) - \mu_t| \geq u) \leq d_{\mathrm{TV}}(\mathcal{L}(X_t), \mathcal{L}(Y)) + \mathbf{Pr}(|f(X_t) - \mu_t| \geq u)$$

$$\leq ne^{-\eta t} + 2e^{-\eta n} + n\exp\left(-\frac{\eta u^2}{nt}\right) + ne^{-\eta nt}.$$

Further, we may assume that $\eta > 0$ is sufficiently small that also Lemma 2.1 holds with this value of $\eta$.

Now let $\kappa = \max\{1, \ln^{-1}(1/\lambda)\}$, let $u \geq 3\kappa\eta^{-1}n^{1/2}\ln n$ and let $t = n^{-1/2}u$. Then

$$\frac{\eta u^2}{nt} = \eta t \geq 3\kappa \ln n.$$

Thus by the above,

$$\mathbf{Pr}(|f(Y) - \mu_t| \geq u) \leq 3ne^{-\eta t} + 2e^{-\eta n}.$$

But $\mu_t = \mathbf{E}[f(X_t)]$, and so by Lemmas 2.1 and 2.4

$$|\mu_t - \mathbf{E}[f(Y)]| \leq d_{\mathrm{W}}(\mathcal{L}(X_t), \mathcal{L}(Y)) = o(1),$$

since $X_0 = \mathbf{0}$ and so $\delta_{n,t} = 0$, and $\eta t \geq 3\kappa \ln n$. Thus we find that

$$\mathbf{Pr}(|f(Y) - \mathbf{E}[f(Y)]| \geq u + 1) \leq 3n\exp(-\eta n^{-1/2}u) + 2e^{-\eta n}$$

$$\leq (3n + 2)\exp(-\eta n^{-1/2}u)$$

for all $n$ sufficiently large and all $3\kappa\eta^{-1}n^{1/2}\ln n \leq u \leq n^{3/2}$.

Let $c = \frac{\eta}{3\kappa}$. Note that $ne^{-cn^{-1/2}u} \geq 1$ for $u \leq 3\kappa\eta^{-1}n^{1/2}\ln n$. Hence, since also $c \leq \eta$,

$$\mathbf{Pr}(|f(Y) - \mathbf{E}[f(Y)]| \geq u + 1) \leq (3n + 2)\exp(-cn^{-1/2}u)$$

for all $n$ sufficiently large and all $0 \leq u \leq n^{3/2}$.

We may assume that $0 < c \leq 1$. Then we may replace $u + 1$ by $u$ if we replace $c$ by $c/2$; that is,

$$\mathbf{Pr}(|f(Y) - \mathbf{E}[f(Y)]| \geq u) \leq (3n + 2)\exp(-(c/2)n^{-1/2}u)$$

for all $n$ sufficiently large and all $0 \le u \le n^{3/2}$. For if $u < 2$, then the right-hand side above is at least 1, and if $u \ge 2$, then $u \ge u/2 + 1$. Now consider square roots: if we replace $c$ by $c/2$ again, then we may replace the factor $(3n+2)$ by $(3n+2)^{1/2}$, which is at most $n$ for $n \ge 4$. Thus, with a new $c$, we have

$$\mathbf{Pr}(|f(Y) - \mathbf{E}[f(Y)]| \ge u) \le n \exp(-cn^{-1/2}u)$$

for all $n$ sufficiently large and all $0 \le u \le n^{3/2}$. One further minor adjustment to $c$ lets us assert that the last inequality holds for all integers $n \ge 2$ and all $0 \le u \le n^{3/2}$.

It remains only to consider values of $u > n^{3/2}$. But always $|f(y)| \le \|y\|_1$, and $\mathbf{E}[\|Y\|_1] \le \frac{\lambda}{1-\lambda}n$. Thus $|f(Y) - \mathbf{E}[f(Y)]| \le \|Y\|_1 + \frac{\lambda}{1-\lambda}n$. Hence for $u > \frac{3\lambda}{1-\lambda}n$,

$$\mathbf{Pr}(|f(Y) - \mathbf{E}[f(Y)]| \ge u) \le \mathbf{Pr}(\|Y\|_1 \ge \mathbf{E}[\|Y\|_1] + u/3) = e^{-\Omega(u)},$$

from the proof of Lemma 2.4 and a standard large deviations calculation for a sum of independent geometric random variables (see [5], pages 201–202), and so Lemma 4.1 follows.

Lemma 4.3 is a quantitative version of some earlier results in [2, 3, 4, 13, 14, 17]. It will also be used in [8] for analyzing the asymptotic distribution of the length of a given queue and for analyzing the "propagation of chaos."

**5. Balance equations and long-term behavior.** In this section we consider the system in equilibrium, and present the key equation (24). We then show that the expected number $\ell(i)$ of queues with at least $i$ customers is close to $n\lambda^{1+d+\cdots+d^{i-1}}$, and that the random number $\ell(i, Y_t)$ of queues with at least $i$ customers stays close to this value over long periods of time. Here $(Y_t)$ is a queue-lengths process in equilibrium. We shall denote $Y_0$ by $Y$ below. As before, $(X_t)$ will denote a queue-lengths process with a "time-dependent" distribution. Observe that $\lambda^{1+d+\cdots+d^{i-1}}$ equals $\lambda^i$ if $d = 1$, and equals $\lambda^{(d^i-1)/(d-1)}$ if $d \ge 2$.

Let $d$ be a fixed positive integer. Fix also a positive integer $n$, and consider the corresponding set $\Omega$ of queue-length vectors. For $x \in \Omega$ and a nonnegative integer $k$, let $u(k, x)$ be the proportion of queues $j$ of length $x(j)$ at least $k$. Thus always $u(0, x) = 1$. Let $u(k)$ denote $\mathbf{E}[u(k, Y)]$: thus $u(k) = \ell(k)/n$. We also let $u_t(k) = \mathbf{E}[u(k, X_t)]$.

If $f$ is the bounded function $f(x) = u(k, x)$, then the generator operator $A$ of the Markov process (see, e.g., [1], Section 1.1, and see also [17]) satisfies

$$Af(x) = \lambda(u(k-1, x)^d - u(k, x)^d) - (u(k, x) - u(k+1, x)).$$

[Compare with (1) earlier.] This is true, since $u(k, x) - u(k+1, x)$ is the proportion of queues of length exactly $k$, and $u(k-1, x)^d - u(k, x)^d$ is the probability that the minimum queue length of the $d$ attempts is exactly $k-1$. From standard theory (see [1], Chapters 1 and 4), for each bounded function $f$, whatever the initial distribution of $X_0$,

$$\frac{d\mathbf{E}[f(X_t)]}{dt} = \mathbf{E}[Af(X_t)],$$

and so in particular, for each positive integer $k$,

(22) $$\frac{du_k(t)}{dt} = \lambda(\mathbf{E}[u(k-1, X_t)^d] - \mathbf{E}[u(k, X_t)^d]) - (u_t(k) - u_t(k+1)).$$

As $Y$ is in equilibrium,

(23) $$\lambda(\mathbf{E}[u(k-1, Y)^d] - \mathbf{E}[u(k, Y)^d]) - (u(k) - u(k+1)) = 0.$$

We shall consider only the equilibrium case for the remainder of this section, and rest our analysis on (23). Now

$$\sum_{k \geq 1} u(k, x) = \frac{1}{n} \sum_{j=1}^{n} x(j) = \frac{1}{n} \|x\|_1,$$

and so

$$\sum_{k \geq 1} u(k) = \frac{1}{n} \mathbf{E}[\|Y\|_1] < \infty.$$

Hence $u(k) \to 0$ as $k \to \infty$. Also $\mathbf{E}[u(k, Y)^d] \leq u(k)$. Summing (23) for $k \geq i$ we obtain, for each $i = 1, 2, \ldots$, that

(24) $$\lambda \mathbf{E}[u(i-1, Y)^d] - u(i) = 0.$$

(This is equivalent to saying that $\mathbf{E}[Af(Y)] = 0$, where $f(x)$ is the number of customers with position at least $i$, but since $f$ is not bounded we cannot assert the result directly.) Equation (24) is crucial to our analysis.

Note first that by (24), $u(1) = \lambda$, and for each $i = 2, 3 \ldots$

$$u(i) = \lambda \mathbf{E}[u(i-1, Y)^d] \leq \lambda u(i-1).$$

Thus by induction on $i$

(25) $$u(i) \leq \lambda^i \qquad \text{for each } i = 0, 1, 2, \ldots.$$

By (24) and the second part of Lemma 4.2, there exists a constant $c_1 > 0$ such that for each positive integer $n$,

(26) $$\sup_i |u(i) - \lambda u(i-1)^d| \leq c_1 n^{-1} \ln^2 n.$$

We claim that, for some constant $c_2 > 0$, for each positive integer $n$

(27) $$\sup_i |u(i) - \lambda^{1+d+\cdots+d^{i-1}}| \leq c_2 n^{-1} \ln^2 n;$$

and hence

(28) $$\sup_i |\ell(i) - n\lambda^{1+d+\cdots+d^{i-1}}| \leq c_2 \ln^2 n$$

for all positive integers $n$.

Let us prove the claim. Note first that for $x, y \geq 0$ we have

(29) $$|y^d - x^d| = |y - x|(y^{d-1} + y^{d-2}x + \cdots + x^{d-1})$$
$$\leq d(x \cup y)^{d-1}|y - x|.$$

Now by (25), $u(i) \leq \lambda^i$ and $\lambda^{1+d+\cdots+d^{i-1}} \leq \lambda^i$, so by (29),

(30) $$|u(i+1) - \lambda^{1+d+\cdots+d^i}|$$
$$\leq |u(i+1) - \lambda u(i)^d| + \lambda |u(i)^d - (\lambda^{1+d+\cdots+d^{i-1}})^d|$$
$$\leq |u(i+1) - \lambda u(i)^d| + \lambda d \cdot \lambda^{i(d-1)}|u(i) - \lambda^{1+d+\cdots+d^{i-1}}|.$$

Since $u(0) = 1$, an easy induction on $i$ using (26) and (30) gives

$$|u(i) - \lambda^{1+d+\cdots+d^{i-1}}| \leq \left(\sum_{j=0}^{i-1}(\lambda d)^j\right) \cdot c_1 n^{-1} \ln^2 n$$

for each $i = 1, 2, \ldots$. Let $i_0$ be the least $i$ such that $\lambda d \cdot \lambda^{i(d-1)} \leq 1/2$. Let $c_2 = c_1 \cdot \max\{2, \sum_{j=0}^{i_0-1}(\lambda d)^j\}$. Clearly

$$|u(i) - \lambda^{1+d+\cdots+d^{i-1}}| \leq c_2 n^{-1} \ln^2 n$$

for $i = 0, 1, \ldots, i_0$. We shall prove by induction on $i$ that this holds for all $i$. Let $i \geq i_0$ and suppose that the inequality holds for $i$. Then by (26) and (30)

$$|u(i+1) - \lambda^{1+d+\cdots+d^i}| \leq c_1 n^{-1} \ln^2 n + \tfrac{1}{2} c_2 n^{-1} \ln^2 n$$
$$\leq c_2 n^{-1} \ln^2 n,$$

as required for the induction step. Thus we have proved the claim (27).

This completes the first half of our task here. Now let $K > 0$ be an arbitrary constant, and let $\tau = n^K$. We see next that all the coordinates $\ell(i, Y_t)$ are likely to stay close to $\ell(i)$ throughout the interval $[0, \tau]$.

LEMMA 5.1. *Let $(Y_t)$ be in equilibrium and let $c > 0$ be a constant. Let $B_\tau$ be the event that for all times $t$ with $0 \leq t \leq \tau$*

$$\sup_i |\ell(i, Y_t) - n\lambda^{1+d+\cdots+d^{i-1}}| \leq cn^{1/2} \ln^2 n.$$

*Then $\mathbf{Pr}(\overline{B_\tau}) \leq e^{-\Omega(\ln^2 n)}$.*

PROOF. By the first part of Lemma 4.2, there exists $\gamma > 0$ such that for all $n$ sufficiently large, for each time $t \geq 0$

$$\mathbf{Pr}\left(\sup_i |\ell(i, Y_t) - \ell(i)| > cn^{1/2} \ln^2 n/4\right) \leq e^{-\gamma \ln^2 n}.$$

Let $s = c\ln^2 n/(8\lambda n^{1/2})$, let $j = \lceil \tau/s \rceil$ and consider times $rs$ for $r = 0, \ldots, j$. The mean number of arrivals in a subinterval $[(r-1)s, rs)$ is $cn^{1/2} \ln^2 n/8$, so by (5) the probability that more than $cn^{1/2} \ln^2 n/4$ arrivals occur is at most $e^{-cn^{1/2} \ln^2 n/24} \leq e^{-\gamma \ln^2 n}$ for $n$ sufficiently large. Then

$$\mathbf{Pr}\left(\sup_i |\ell(i, Y_t) - \ell(i)| > cn^{1/2} \ln^2 n/2 \text{ for some } t \in [0, \tau]\right)$$

$$\leq \sum_{r=0}^{j} \mathbf{Pr}\left(\sup_i |\ell(i, Y_{rs}) - \ell(i)| \geq cn^{1/2} \ln^2 n/4\right)$$

$$+ \sum_{r=1}^{j} \mathbf{Pr}([(r-1)s, rs) \text{ has } > cn^{1/2} \ln^2 n/4 \text{ arrivals})$$

$$\leq (\tau/s + 2) \cdot 2e^{-\gamma \ln^2 n}$$

$$\leq e^{-(\gamma/2) \ln^2 n}$$

for $n$ sufficiently large. Now we can use (28). □

**6. Lower bound on mixing times.** In this short section we show that the upper bounds on the mixing times in Theorem 1.1 are of the right order, and in particular we prove (2). We do this by considering the total number of nonempty queues in the system. The idea is that this number is highly concentrated around its mean. In equilibrium the mean is $\lambda n$. On the other hand, if $X_0 = \mathbf{0}$ and $t \leq \theta \ln n$, then the expected number of nonempty queues is less than $(\lambda - \theta)n$ for a suitable constant $\theta$. Thus for such $t$ the two distributions are far apart.

Now for the details. Consider two $n$-queue processes $(X_t)$ and $(Y_t)$, where $X_0 = \mathbf{0}$ and $(Y_t)$ is in equilibrium. By Lemma 2.3, we can couple $(X_t)$ and $(Y_t)$ in such a way that always $X_t \leq Y_t$ and so $u(i, X_t) \leq u(i, Y_t)$ for each $i$. In particular, if as before we let $u_t(i) = \mathbf{E}[u(i, X_t)]$, then $u_t(i) \leq u(i)$. (Recall that $u(i) = \mathbf{E}[u(i, Y_t)] = \frac{1}{n}\mathbf{E}[\ell(i, Y_t)]$.) Let $s_t(i) = u(i) - u_t(i)$, so $s_t(i) \geq 0$. Also from (22)

$$\frac{du_t(1)}{dt} = \lambda(1 - \mathbf{E}[u(1, X_t)^d]) - (u_t(1) - u_t(2))$$

and

$$0 = \lambda(1 - \mathbf{E}[u(1, Y_t)^d]) - (u(1) - u(2)),$$

so that
$$\frac{ds_t(1)}{dt} = -s_t(1) - \lambda(\mathbf{E}[u(1,Y_t)^d] - \mathbf{E}[u(1,X_t)^d]) + s_t(2).$$

But, as in (29), given $0 \leq y \leq x \leq 1$, we have
$$x^d - y^d = (x-y)(x^{d-1} + \cdots + y^{d-1}) \leq d(x-y).$$

Hence, since we may assume that always $u(1,X_t) \leq u(1,Y_t)$, we have
$$\mathbf{E}[u(1,Y_t)^d - u(1,X_t)^d] \leq d(u(1) - u_t(1)) = ds_t(1).$$

It follows [using also the fact that $s_t(2) \geq 0$] that for all $t \geq 0$,
$$\frac{ds_t(1)}{dt} \geq -(1+\lambda d)s_t(1).$$

Thus
$$s_t(1) \geq \lambda e^{-(1+\lambda d)t}$$
for all times $t \geq 0$, since $s_0(1) = u(1) = \lambda$. Thus if $t \leq \frac{1}{2(1+\lambda d)} \ln n - \frac{2}{1+\lambda d} \ln \ln n$, then $s_1(t) \geq \lambda n^{-1/2} \ln^2 n$, that is,
$$\lambda n - \mathbf{E}[\ell(1,X_t)] \geq \lambda n^{1/2} \ln^2 n.$$

But if $t$ is $O(\ln n)$, then from (21) in the proof of Lemma 4.3,
$$\mathbf{Pr}(|\ell(1,X_t) - \mathbf{E}[\ell(1,X_t)]| \geq n^{1/2} \ln^{3/2} n) = e^{-\Omega(\ln^2 n)}.$$

Also, by the first part of Lemma 4.2,
$$\mathbf{Pr}(|\ell(1,Y_t) - \lambda n| \geq \tfrac{1}{2}\lambda n^{1/2} \ln^2 n) = e^{-\Omega(\ln^2 n)}.$$

Inequality (2) now follows.

**7. Proof of Theorem 1.3.** We assume throughout that the process is in equilibrium. For each time $t \geq 0$ and each $i = 0, 1, \ldots$ let the random variable $Z_t(i)$ be the number of new customers arriving during $[0,t]$ with position at least $i$ on arrival. Let $J_0 = 0$ and enumerate the arrival times after time 0 as $J_1, J_2, \ldots$. We define a "horizon" time $t_0 = \ln^2 n$, and let $N = \lceil 2\lambda n t_0 \rceil$. Note that $\ell(1,Y_t)$ is the number of nonempty queues at time $t$. For each time $t$, let $A_t$ be the event
$$\{\lambda n/2 \leq \ell(1,Y_s) \leq 2\lambda n \ \forall\, s \in [0,t]\}.$$

Then by Lemma 5.1, the event $A_{t_0}$ holds a.a.s.

We need two more preliminary results. The first one is an analogue of a special case of Lemma 2.1 in [7]. It may be proved quickly along the lines of the proof of that lemma; for completeness we give a proof here.

LEMMA 7.1. *Let $(X_t)$ be in equilibrium. Let $s, \tau > 0$ and let $a, b$ be non-negative integers. Let $\delta = n(\lambda ds)^{b+1}/(b+1)!$. Then*

$$(31) \quad \mathbf{Pr}[M_t \leq a \text{ for some } t \in [0, \tau]] \leq \left(\frac{\tau}{s} + 1\right)(\mathbf{Pr}(M_0 \leq a + b) + \delta)$$

*and*

$$(32) \quad \mathbf{Pr}[M_t \geq a + b \text{ for some } t \in [0, \tau]] \leq \left(\frac{\tau}{s} + 1\right)(\mathbf{Pr}(M_0 \geq a) + \delta).$$

PROOF. The $j = \lfloor \frac{\tau}{s} \rfloor + 1$ disjoint intervals $[(r-1)s, rs)$ for $r = 1, \ldots, j$ cover $[0, \tau]$. Let $C_r$ denote the event of having at least $b+1$ arrivals in the interval $[(r-1)s, rs)$ which are placed into a single bin. Then

$$\mathbf{Pr}(C_r) \leq n\mathbf{Pr}(Po(\lambda ds) \geq b + 1) \leq \delta.$$

But

$$\{M_t \leq a \text{ for some } t \in [0, \tau]\} \subseteq \left(\bigcup_{r=1}^{j}\{M_{rs} \leq a + b\}\right) \cup \left(\bigcup_{r=1}^{j} C_r\right)$$

and (31) follows. Similarly

$$\{M_t \geq a + b \text{ for some } t \in [0, \tau)\} \subseteq \left(\bigcup_{r=0}^{j-1}\{M_{rs} \geq a\}\right) \cup \left(\bigcup_{r=1}^{j} C_r\right)$$

and (32) follows. □

We shall use the second lemma to bound "initial" effects. For each time $t \geq 0$ let $S_t$ be the event that some initial customer survives at least to time $t$.

LEMMA 7.2. *Let $0 < \lambda < 1$ and let $d$ be a positive integer. Let $\alpha = \min\{\frac{1}{4}\ln(\frac{1}{\lambda}), \frac{1}{4}\}$, so $\alpha > 0$. Then for each positive integer $n$ and each time $t \geq 0$,*

$$\mathbf{Pr}(S_t) \leq 2ne^{-\alpha t}.$$

PROOF. Let $k = \lfloor t/4 \rfloor$. By (25)

$$\mathbf{Pr}(M_0 \geq k + 1) \leq n\lambda^{k+1} \leq n\lambda^{t/4}.$$

Let $Z$ be distributed like the sum of $k$ independent service times. Thus $Z$ has a gamma distribution, and has moment generating function $\mathbf{E}[e^{sZ}] =$

$(1-s)^{-k}$ for $s < 1$. So by Markov's inequality, for each $0 \leq s < 1$,
$$\mathbf{Pr}(Z \geq t) \leq e^{-ts}(1-s)^{-k}$$
$$= e^{-t/2} 2^k \qquad \text{taking } s = \frac{1}{2}$$
$$\leq \left(\frac{e^2}{2}\right)^{-t/4} \leq e^{-t/4}.$$

Hence
$$\mathbf{Pr}(S_t) \leq \mathbf{Pr}(M_0 \geq k+1) + n\mathbf{Pr}(Z \geq t) \leq n\lambda^{t/4} + ne^{-t/4} \leq 2ne^{-\alpha t}. \quad \square$$

Let $i^* = i^*(n)$ be the smallest integer such that $\lambda^{(d^i-1)/(d-1)} < n^{-1/2}\ln^2 n$ (see the discussion following Theorem 1.3). Note that $i^* = \ln\ln n / \ln d + O(1)$.

We can handle the lower bound on $M_t$ easily. By Lemma 5.1, $\mathbf{Pr}(M < i^* - 1) = e^{-\Omega(\ln^2 n)}$. In particular $M \geq i^* - 1$ a.a.s. Further, (31), with $\tau = n^K$, $a = i^* - 3$, $b = 1$ and $s = n^{-K-2}$, shows that

(33) $$\min\{M_t : 0 \leq t \leq n^K\} \geq i^* - 2 \qquad \text{a.a.s.,}$$

since
$$\delta = n\mathbf{Pr}(Po(\lambda ds) \geq 2) \leq n(\lambda ds)^2 = O(n^{-2K-3}).$$

This result establishes the lower bound half of (3).

The upper bounds on $M_t$ are less straightforward to prove. We consider first the easier case when $d \geq 3$.

7.1. *The case $d \geq 3$.* We shall show that $M \leq i^*$ a.a.s. For $k = 0, 1, \ldots$, let $E_k$ be the event that $\ell(i^*, Y_{J_k}) \leq 2n^{1/2}\ln^2 n$; that is, at time (just after) $J_k$ there are no more than $2n^{1/2}\ln^2 n$ queues with at least $i^*$ customers. Then $\mathbf{Pr}(\overline{E_k}) = e^{-\Omega(\ln^2 n)}$ by Lemma 5.1. Consider the customer who arrives at time $J_k$: on $E_{k-1}$, he has probability at most $p_1 = (2n^{-1/2}\ln^2 n)^d$ of joining a queue of length at least $i^*$. Note that
$$\mathbf{Pr}(J_{N+1} \leq t_0) \leq \mathbf{Pr}[Po(\lambda n t_0) \geq 2\lambda n t_0] = e^{-\Omega(n \ln^2 n)}.$$

Also, for each positive integer $r$,
$$\mathbf{Pr}(B(N, p_1) \geq r) \leq (Np_1)^r = O((n^{-(d/2-1)}(\ln n)^{2d+2})^r).$$

Hence, for each positive integer $r$,
$$\mathbf{Pr}[Z_{t_0}(i^* + 1) \geq r]$$
$$\leq \mathbf{Pr}(B(N, p_1) \geq r) + \mathbf{Pr}\left(\bigcup_{k=0}^{N-1} \overline{E_k}\right) + \mathbf{Pr}(J_{N+1} \leq t_0)$$
$$= O((n^{-(d/2-1)}(\ln n)^{2d+2})^r).$$

Also, by Lemma 7.2 there exists a constant $\alpha > 0$ such that the probability that some "initial" customer has not departed by time $t_0$ is at most $2ne^{-\alpha t_0}$. Hence, for each positive integer $r$,

$$\mathbf{Pr}[M \geq i^* + r] \leq \mathbf{Pr}[Z_{t_0}(i^* + 1) \geq r] + 2ne^{-\alpha t_0}.$$

In particular, $\mathbf{Pr}(M \geq i^* + 1) = o(1)$, which together with the earlier result that $M \geq i^* - 1$ a.a.s. completes the proof that $M$ is concentrated on at most the two values $i^*$ and $i^* - 1$. Also,

(34) $$\mathbf{Pr}[M \geq i^* + 2K + 5] = o(n^{-K-2}).$$

Now (32), with $\tau = n^K$, $a = i^* + 2K + 5$, $b = \lceil K/2 \rceil + 1$ and $s = n^{-2}$, lets us complete the proof of (3), since

$$\delta = n\mathbf{Pr}(Po(\lambda ds) \geq b + 1) \leq n(\lambda ds)^{K/2+2} = O(n^{-K-3}).$$

7.2. *The case $d = 2$.* The case $d = 2$ needs a little more effort, and uses the "drift" results from Section 7 in [7]. For convenience we restate these results here, as the next two lemmas. The first lemma concerns hitting times for a generalized random walk with "drift."

LEMMA 7.3. *Let $\phi_0 \subseteq \phi_1 \subseteq \cdots \subseteq \phi_m$ be a filtration, and let $Y_1, Y_2, \ldots, Y_m$ be random variables taking values in $\{-1, 0, 1\}$ such that each $Y_i$ is $\phi_i$-measurable. Let $E_0, E_1, \ldots, E_{m-1}$ be events where $E_i \in \phi_i$ for each $i$, and let $E = \bigcap_i E_i$. Let $R_t = R_0 + \sum_{i=1}^{t} Y_i$. Let $0 \leq p \leq 1/3$, let $r_0$ and $r_1$ be integers such that $r_0 < r_1$ and let $pm \geq 2(r_1 - r_0)$. Assume that for each $i = 1, \ldots, m$,*

$$\mathbf{Pr}(Y_i = 1 | \phi_{i-1}) \geq 2p \quad \text{on } E_{i-1} \cap \{R_{i-1} < r_1\}$$

*and*

$$\mathbf{Pr}(Y_i = -1 | \phi_{i-1}) \leq p \quad \text{on } E_{i-1} \cap \{R_{i-1} < r_1\}.$$

*Then*

$$\mathbf{Pr}(E \cap \{R_t < r_1 \ \forall t \in \{1, \ldots, m\}\} | R_0 = r_0) \leq \exp\left(-\frac{pm}{28}\right).$$

The second lemma shows that if we try to cross an interval against the drift, then we will rarely succeed.

LEMMA 7.4. *Let $a$ be a positive integer. Let $p$ and $q$ be reals with $q > p \geq 0$ and $p + q \leq 1$. Let $\phi_0 \subseteq \phi_1 \subseteq \phi_2 \subseteq \cdots$ be a filtration, and let $Y_1, Y_2, \ldots$ be random variables taking values in $\{-1, 0, 1\}$ such that each $Y_i$ is $\phi_i$-measurable.*

Let $E_0, E_1, \ldots$ be events where each $E_i \in \phi_i$, and let $E = \bigcap_i E_i$. Let $R_0 = 0$ and let $R_k = \sum_{i=1}^{k} Y_i$ for $k = 1, 2, \ldots$. Assume that for each $i = 1, 2, \ldots$,

$$\mathbf{Pr}(Y_i = 1 | \phi_{i-1}) \leq p \qquad \text{on } E_{i-1} \cap \{0 \leq R_{i-1} \leq a - 1\}$$

and

$$\mathbf{Pr}(Y_i = -1 | \phi_{i-1}) \geq q \qquad \text{on } E_{i-1} \cap \{0 \leq R_{i-1} \leq a - 1\}.$$

Let

$$T = \inf\{k \geq 1 : R_k \in \{-1, a\}\}.$$

Then

$$\mathbf{Pr}(E \cap \{R_T = a\}) \leq (p/q)^a.$$

Having stated the two "drift" lemmas, let us resume the proof for the case $d = 2$. Still let $t_0 = \ln^2 n$ and $N = \lceil 2\lambda n t_0 \rceil$. We first show that $M \geq i^*$, by showing that in fact

(35) $$\ell(i^*, Y_{t_0}) \geq \ln^3 n \qquad \text{a.a.s.}$$

Let $J'_0 = 0$, and enumerate all jump times after time 0 (not just the arrival times) as $J'_1, J'_2, \ldots$. For $k = 0, 1, \ldots$ let $E_k$ be the event that $\ell(i^* - 1, Y_{J'_k}) \geq \frac{1}{2} n^{1/2} \ln^2 n$; and let $E'_k$ be the event that $\ell(i^*, Y_{J'_k}) \leq 2 \ln^3 n$. By Lemma 5.1 as before, $\mathbf{Pr}(\overline{E_k}) = e^{-\Omega(\ln^2 n)}$. Let $V_k = \ell(i^*, Y_{J'_k}) - \ell(i^*, Y_{J'_{k-1}})$ for $k = 1, 2, \ldots$, so that

$$\ell(i^*, Y_{J'_k}) = \ell(i^*, Y_0) + \sum_{j=1}^{k} V_j.$$

Let $p_2 = \ln^4 n / (25n)$. On $A_{J'_{k-1}} \cap E_{k-1} \cap E'_{k-1}$,

$$\mathbf{Pr}(V_k = 1 | \phi_{J'_{k-1}}) \geq 2p_2$$

and

$$\mathbf{Pr}(V_k = -1 | \phi_{J'_{k-1}}) \leq p_2,$$

for $n$ sufficiently large.

Also then $np_2 \geq 4 \ln^3 n$ for $n$ sufficiently large. Hence by Lemma 7.3, a.a.s. $\ell(i^*, Y_{J'_i}) \geq 2 \ln^3 n$ for some $i \leq n$. Also $J'_n \leq t_0$ a.a.s., since

$$\mathbf{Pr}(J'_n > t_0) \leq \mathbf{Pr}(J_n > t_0) = \mathbf{Pr}(Po(\lambda n t_0) < n) = e^{-\Omega(n \ln^2 n)}.$$

Thus a.a.s. $\ell(i^*, Y_t) \geq 2 \ln^3 n$ for some $t \in [0, t_0]$.

It now suffices to show that a.a.s. there will be no "excursions" that cross downward from $\lceil 2 \ln^3 n \rceil$ to at most $\ln^3 n$ during the interval $[0, \tau]$. Let $B$ be the event that there is such a crossing. The only possible start times for such a crossing are departure times during $[0, t_0]$, and a.a.s. there are at most $N$ such times. We may use Lemma 7.4 to upper bound the probability that any given excursion leads to a crossing. Let $a = \lceil \ln^3 n \rceil$. Let $p = p_2$ and let $q = 2p_2$. We apply the lemma with $Y_k$ replaced by $-V_k$ and with $a$, $p$ and $q$ as above, and with the events $\tilde{E}_k = A_{J'_k} \cap E_k$. We obtain

$$\mathbf{Pr}(B) \leq N 2^{-a} + \mathbf{Pr}\left(\bigcup_{k=0}^{N-1} \overline{\tilde{E}_k}\right) + o(1) = o(1).$$

Thus we have proved that $M \geq i^*$ a.a.s.

We now consider upper bounds on $M$. We shall show that $M \leq i^* + 1$ a.a.s., by showing that $\ell(i^* + 2, Y_{t_0}) = 0$ a.a.s. For $k = 0, 1, \ldots$, let $F_k$ be the event that at the arrival time $J_k$ there are no more than $2n^{1/2} \ln^2 n$ queues of length at least $i^*$. By Lemma 5.1 we have $\mathbf{Pr}(\overline{F_k}) = e^{-\Omega(\ln^2 n)}$. Consider the customer arriving at time $J_k$: on $F_{k-1}$ he has probability at most $p_3 = 4 \ln^4 n / n$ of joining a queue of length at least $i^*$. Thus for each positive integer $r$,

$$\mathbf{Pr}(Z_{t_0}(i^* + 1) \geq r) \leq \mathbf{Pr}(B(N, p_3) \geq r) + \mathbf{Pr}\left(\bigcup_{k=0}^{N-1} \overline{F_k}\right) + \mathbf{Pr}(J_{N+1} \leq t_0).$$

Also, by Lemma 7.2 the probability that some "initial" customer has not departed by time $t_0/2$ is at most $2n e^{-\alpha t_0/2}$ for a suitable constant $\alpha > 0$. Hence, there is a constant $\tilde{c}$ such that with probability $1 - e^{-\Omega(\ln^2 n)}$ we have $\ell(i^* + 1, Y_t) \leq \tilde{c} \ln^6 n$ uniformly for all $t \in [t_0/2, t_0]$. Thus this also holds over $[0, t_0]$.

For $k = 0, 1, \ldots$, let $F'_k$ be the event that $\ell(i^* + 1, Y_{J_k}) \leq \tilde{c} \ln^6 n$; that is, at time $J_k$ there are no more than $\tilde{c} \ln^6 n$ queues with at least $i^* + 1$ customers. On $F'_{k-1}$, the customer arriving at time $J_k$ has probability at most $p_4 = \tilde{c}^2 (\ln n)^{12} n^{-2}$ of joining a queue of length at least $i^* + 1$. Then for each positive integer $r$,

$$\mathbf{Pr}(Z_{t_0}(i^* + 2) \geq r) \leq \mathbf{Pr}(B(N, p_4) \geq r) + \mathbf{Pr}\left(\bigcup_{k=0}^{N-1} \overline{F'_k}\right) + \mathbf{Pr}(J_{N+1} \leq t_0).$$

Also, by Lemma 7.2, the probability that some initial customer stays until time $t_0$ is at most $2n e^{-\alpha t_0}$. It follows that a.a.s. $M_{t_0} \leq i^* + 1$, and

(36) $$\mathbf{Pr}(M_{t_0} \geq i^* + K + 5) = O(n^{-K-3}).$$

Now the inequality (32) in Lemma 7.1 with $\tau = n^K$, $a = i^* + K + 5$, $b = \lceil K/2 \rceil + 1$ and $s = n^{-2}$ yields the upper bound part of (3), since

$$\delta = n \mathbf{Pr}(Po(\lambda ds) \geq b + 1) \leq n(\lambda ds)^{b+1} = O(n^{-K-3}).$$

Finally, let us note one result which will be useful in [8]. Let the integer $d \geq 2$ be fixed. Then, as in (36), if $r = O(\ln n)$,

$$\mathbf{Pr}(M \geq i^* + 1 + r) = e^{-\Omega(r \ln n)}. \tag{37}$$

**8. Concluding remarks.** We have investigated the well-known "supermarket" model with $n$ servers and a fixed number $d$ of random choices. We have shown that the system converges rapidly to its equilibrium distribution. Our main result is that, for $d \geq 2$, in equilibrium the maximum length of a queue is a.a.s. concentrated on two adjacent values close to $\ln \ln n / \ln d$. In contrast, when $d = 1$ the maximum length of a queue is a.a.s. close to $\ln n / \ln \ln n$. Along the way we used the fact (27) that, in equilibrium, for each nonnegative integer $k$ the proportion of queues of length at least $k$ is close to $\lambda^{1+d+\cdots+d^{k-1}}$.

In [8] we use this last result together with mixing and concentration estimates (and upper bounds on the maximum length of a queue) obtained here to prove quantitative results on the convergence of the distribution of a queue length. In particular, we give a quantitative version of the convergence result mentioned in the Introduction. Let $(v_t(k): k \in \mathbb{N})$ be as in (1) above. It turns out that, for suitable initial conditions, uniformly over all $t \geq 0$ the distribution of a given queue length at time $t$ is close in total variation distance to the distribution of an integer-valued nonnegative random variable $V_t$ such that $\mathbf{Pr}(V_t \geq k) = v_t(k)$. Also in [8] we investigate the asymptotic independence of small subsets of queues, the "chaotic behavior" of the system.

DEPARTMENT OF MATHEMATICS
LONDON SCHOOL OF ECONOMICS
HOUGHTON STREET
LONDON WC2A 2AE
UNITED KINGDOM
E-MAIL: m.j.luczak@lse.ac.uk

DEPARTMENT OF STATISTICS
UNIVERSITY OF OXFORD
1 SOUTH PARKS ROAD
OXFORD OX1 3TG
UNITED KINGDOM
E-MAIL: cmcd@stats.ox.ac.uk